\documentclass[final,3p]{elsarticle}
\usepackage[cmex10]{amsmath}
\usepackage[bitstream-charter]{mathdesign}
\usepackage[T1]{fontenc}

\usepackage{lineno}
\usepackage{fancyhdr}
\usepackage{fixltx2e}

\usepackage{array}
\usepackage{mdwtab}
\usepackage{eqparbox}
\usepackage{url}
\usepackage{color}
\usepackage{multirow}
\usepackage{tikz,pgfplots}
\usepackage[%
  plainpages=false,
  colorlinks=true,
  citecolor=blue,
  filecolor=blue,
  linkcolor=blue,
  urlcolor=blue
  ]{hyperref}

\newcommand{\note}[1]{\noindent\emph{\textcolor{blue}{#1}}}

\newcommand{\smallsub}[1]{\!\mbox{\scriptsize #1}}

\usepackage[shortcuts]{extdash}
\usepackage{booktabs}
\usepackage[mathscr]{euscript}
\usepackage{color}
\usepackage{graphicx}
\usepackage{algorithmic,algorithm}
\usepackage{tikz,pgfplots}
\usepackage{upgreek}
\usepackage{multirow}
\usepackage{yfonts}
\usepackage{mathtools}
\usepackage[normalem]{ulem}
\usepackage{latexsym}
\usepackage{url}

\newcommand{\zapspace}{\topsep=0pt\partopsep=0pt\itemsep=0pt\parskip=0pt}

\definecolor{utorange}{rgb}{0.8,0.33,0.}
\definecolor{themec}{RGB}{51,108,121}
\definecolor{darkred}{rgb}{.6,.1,.1}
\definecolor{darkblue}{rgb}{.1,.1,.9}
\definecolor{greenback}{rgb}{.19,.94,.13}
\definecolor{orange}{rgb}{.76,.39,.13}
\definecolor{grass}{rgb}{.19,.64,.13}
\definecolor{sierp}{RGB}{209,28,209}
\definecolor{bgorange}{rgb}{1.,.95,.78}
\definecolor{grassgreen}{RGB}{92,135,39}
\definecolor{thinbox}{rgb}{.7,.8,1.}

\renewcommand{\vec}[1]{{\mathchoice
                     {\mbox{\boldmath$\displaystyle{#1}$}}
                     {\mbox{\boldmath$\textstyle{#1}$}}
                     {\mbox{\boldmath$\scriptstyle{#1}$}}
                     {\mbox{\boldmath$\scriptscriptstyle{#1}$}}}}

\newcommand\restr[2]{{
  \left.\kern-\nulldelimiterspace 
  {#1}\vphantom{\big|} \right|_{#2}}}

\newcommand{\B}{\mathcal{B}}

\newcommand{\J}{\mathcal{J}}

\newcommand{\Reg}{\mathcal{R}}

\newcommand{\bit}{\begin{itemize}}
\newcommand{\eit}{\end{itemize}}
\newcommand{\bdm} {\begin{displaymath}}
\newcommand{\edm} {\end{displaymath}}

\newcommand{\Cprior}{\mathcal{C}_{\text{\tiny{prior}}}}

\newcommand{\observables}{\bs d}
\newcommand{\data}{ \observables_{\rm{obs}} }
\newcommand{\uobs}{ \uu_{\rm{obs}} }

\newcommand{\FF}{{\ensuremath{\matrix{F}}}}
\newcommand{\W}{{\ensuremath{\matrix{W}}}}

\newcommand{\ncov} {\bs{\Gamma}_{\!\mbox{\tiny noise}} }
\newcommand{\prcov} {\bs{\Gamma}_{\!\mbox{\tiny prior}} }
\newcommand{\postcov} {\bs{\Gamma}_{\!\mbox{\tiny post}} }
\newcommand{\map} {{\beta}_{\mbox{\tiny MAP}} }
\newcommand{\dmap} {{\vec{\beta}}_{\mbox{\tiny MAP}} }
\newcommand{\mpr} {{\vec{\beta}}_{\mbox{\tiny prior}} }

\newcommand{\Hmisfit}{ \matrix{H}_{\mbox{\tiny misfit}} }
\newcommand{\HT}{\matrix{\tilde{H}}_{\mbox{\tiny misfit}}}

\newcommand{\SI}{ \ensuremath{{\Sigma}} }
\newcommand{\PP}{ \ensuremath{\matrix{\Psi}} }
\newcommand{\WW}{ \mathcal{W}}

\newcommand{\mc}[1]{\mathcal{#1}}

\newcommand{\gbf}[1]{\boldsymbol{#1}}
\newcommand{\bs}[1]{\ensuremath{\boldsymbol{#1}}}

\newcommand{\edot}{\dot{\gbf{\varepsilon}}}

\newcommand{\secinve}{\edot_\mathrm{II}}

\DeclareMathOperator{\diag}{diag}

\newcommand{\uh}{\hat{\uu}}
\newcommand{\ph}{\hat{p}}
\newcommand{\vh}{\hat{\vv}}
\newcommand{\qh}{\hat{q}}
\newcommand{\betah}{\hat{\beta}}

\renewcommand{\vec}[1]{\gbf{#1}}

\newcommand{\uu}{\ensuremath{\vec{u}}}
\newcommand{\vv}{\ensuremath{\vec{v}}}

\renewcommand{\matrix}[1] {\ensuremath{\boldsymbol{#1}}}

\modulolinenumbers[5]

\journal{Journal of Computational Physics}

\begin{document}

\begin{frontmatter}

\title{Scalable and efficient algorithms for the propagation of
  uncertainty \\ from data through inference to prediction for
  large-scale problems,\\ with
  application to flow of the Antarctic ice sheet\\[-1.5ex]}
\author[addressices]{Tobin Isaac}
\ead{tisaac@ices.utexas.edu}
\author[addressucm]{Noemi Petra
\corref{mycorrespondingauthor}}
\cortext[mycorrespondingauthor]{Corresponding author}
\ead{npetra@ucmerced.edu}
\author[addressnyu]{Georg Stadler}
\ead{stadler@cims.nyu.edu}
\author[addressices,addressmech,addressgeo]{Omar Ghattas}
\ead{omar@ices.utexas.edu}

\address[addressices]{Institute for Computational Engineering \&
  Sciences, The University of Texas at Austin, Austin, TX, USA}
\address[addressucm]{Applied Mathematics, School of Natural Sciences, University of California,
  Merced, CA, USA}
\address[addressnyu]{Courant Institute of Mathematical Sciences, New
  York University, New York, NY, USA}
\address[addressmech]{Department of Mechanical Engineering, The
  University of Texas at Austin, Austin, TX, USA}
\address[addressgeo]{Jackson School of Geosciences, The
  University of Texas at Austin, Austin, TX, USA}

\begin{abstract}
The majority of research on efficient and scalable algorithms in
computational science and engineering has focused on the {\em forward
  problem}: given parameter inputs, solve the governing equations to
determine output quantities of interest. In contrast, here we consider
the broader question: given a (large-scale) model containing uncertain
parameters, (possibly) noisy observational data, and a prediction
quantity of interest, how do we construct efficient and scalable
algorithms to (1) infer the model parameters from the data (the
{\em deterministic inverse problem}), (2) quantify the uncertainty in the
inferred parameters (the {\em Bayesian inference problem}), and (3)
propagate the resulting uncertain parameters through the model to
issue predictions with quantified uncertainties (the {\em forward
uncertainty propagation problem})?

We present efficient and scalable algorithms for this end-to-end,
data-to-prediction process under the Gaussian approximation and in the
context of modeling the flow of the Antarctic ice sheet and its
effect on loss of grounded ice to the ocean. The ice is modeled as a viscous, incompressible,
creeping, shear-thinning fluid. The observational data come from %
satellite measurements of surface ice flow velocity, and the uncertain
parameter field to be inferred is the basal sliding parameter,
represented by a heterogeneous coefficient in a Robin boundary
condition at the base of the ice sheet. The prediction quantity of
interest is the present-day ice mass flux from the Antarctic continent to the
ocean.

We show that the work required for executing this data-to-prediction
process---measured in number of forward (and adjoint) ice sheet model
solves---is independent of the state dimension, parameter dimension,
data dimension, and the number of processor cores. The key to achieving
this dimension independence is to exploit the fact that, despite their
large size, the observational data typically provide only sparse
information on model parameters. This property can be exploited to
construct a low rank approximation of the linearized
parameter-to-observable map via randomized SVD methods and
adjoint-based actions of Hessians of the data misfit functional.

\end{abstract}


\end{frontmatter}

\section{Introduction}
\label{sec:intro}

The future mass balance of the polar ice sheets will be critical to
climate in the coming century, yet there is much uncertainty surrounding
even their current mass balance.  The current rate of ice sheet mass
loss was recently estimated at roughly 200 billion metric tons per year
in \cite{ShepherdIvinsGeruoEtAl12} using data from various sources,
including radar and laser altimetry, gravimetric observations, and
surface mass balance calculations of regional climate models.\footnote{This
estimate is broken down into estimates for individual ice sheets with
confidence intervals ($-149 \pm 49$ Gigatonnes per year from Greenland,
$+14 \pm 43$ from East Antarctica, $-65 \pm 26$ from West Antarctica,
$-20 \pm 15$ from the Antarctic peninsula).} Moreover, this mass
loss has been observed to be accelerating \cite{HannaNavarroPattynEtAl13}.  Driven by increased warming, collapse
of even a small portion of one of these ice sheets has the potential to
greatly accelerate this figure.  Indeed, recent evidence suggests that
sea level rose abruptly at the end of the last interglacial period
(118,000 years ago) by 5--6 m; the likely cause is catastrophic collapse
of an ice sheet driven by warming oceans
\cite{OlearyHeartyThompsonEtAl13}. Based on a conservative estimate of a
half meter of sea level rise, the Organization for Economic Cooperation
and Development estimates that the 136 largest port cities, with 150
million inhabitants and \$35 trillion worth of assets, will be at risk
from coastal flooding by 2070 \cite{NichollsHansonHerweijerEtAl08}.

Clearly, model-based projections of the evolution of the polar ice
sheets will play a central role in anticipating future sea level
rise. However, current ice sheet models are subject to considerable
uncertainties. Indeed, ice sheet models were left out of the
Intergovernmental Panel on Climate Change's 4th Assessment Report,
which stated that ``the uncertainty in the projections of the land ice
contributions [to sea level rise] is dominated by the various
uncertainties in the land ice models themselves...rather than in the
temperature projections'' \cite{IPCC07a}.

Ice is modeled as a creeping, viscous,
incompressible, shear-thinning fluid with strain-rate-
and temperature-dependent viscosity. Severe mathematical and
computational challenges
place significant barriers on improving predictability of ice sheet
flow models. These include complex and very high-aspect ratio (thin)
geometry, highly nonlinear and anisotropic rheology, extremely
ill-conditioned linear and nonlinear algebraic systems that arise upon
discretization as a result of heterogeneous, widely-varying viscosity and basal
sliding parameters, a broad range of relevant length scales (tens of meters to thousands
of kilometers), localization phenomena including fracture, and complex
sub-basal hydrological processes.

However, the greatest mathematical and computational challenges lie in
quantifying the uncertainties in the predictions of the ice sheet
models. These models are characterized by unknown or poorly
constrained  fields
describing the basal sliding parameter (resistance to sliding at the base of
the ice sheet), basal topography, geothermal heat flux, and
rheology. 
While many of these parameter fields cannot be directly observed, they can be
inferred from satellite observations, such as those of ice surface velocities or ice thickness, which
leads to a severely ill-posed inverse problem whose solution is
extremely challenging. Quantifying the uncertainties that result from
inference of these ice sheet parameter fields from noisy data can be
accomplished via the framework of Bayesian inference. Upon
discretization of the unknown infinite dimensional parameter field,
the solution of the Bayesian inference problem takes the form of a
very high-dimensional {\em posterior} probability density function
(pdf) that assigns to any candidate set of parameter fields our belief
(expressed as a probability) that a member of this candidate set is
the ``true'' parameter field that gave rise to the observed
data. Sampling this posterior pdf to compute, for example, the mean
and covariance of the parameters presents tremendous challenges, since
not only is it high-dimensional, but evaluating the pdf at any point
in parameter space requires a forward ice sheet flow simulation---and
millions of such evaluations may be required to obtain statistics of
interest using state-of-the-art Markov chain Monte Carlo
methods. Finally, the ice sheet model parameters and their associated
uncertainties can be propagated through the ice sheet flow model to
yield predictions of not only the mass flux of ice into the ocean, but
also the confidence we have in those predictions. This amounts to
solving a system of stochastic PDEs, which again is intractable when
the PDEs are complex and highly nonlinear and the parameters are
high-dimensional due to discretization of an infinite-dimensional
field.

In summary, while one can {\em formulate} a data-to-prediction
framework to quantify uncertainties from data to inferred model
parameters to predictions with an underlying model of non-Newtonian
ice sheet flow, attempting to {\em execute} this framework for the
Antarctic ice sheet (or other large-scale complex models) is
intractable for high-dimensional parameter fields using current
algorithms.
Yet, quantifying the uncertainties in predictions of ice sheet models
is essential if these models are to play a significant role in
projections of future sea level.  The purpose of this paper is to
present an integrated framework and efficient, scalable algorithms for
carrying out this data-to-prediction process. By {\em scalable}, we
mean that the cost---measured in number of (linearized) forward (and
adjoint) solves---is independent of
not only the number of processor cores, but importantly the state
variable dimension, the parameter dimension, and the data dimension.

Two key ideas are needed to produce such scalable algorithms. First, we
use Gaussian approximations of both the posterior pdf that results from
Bayesian solution of the inverse problem of inferring ice sheet
parameter fields from satellite observations of surface velocity, as
well as the pdf resulting from propagating the uncertain parameter
fields through the forward ice sheet model to yield predictions of
present-day mass flux into the ocean. %
This is
accomplished by linearizing the parameter-to-observable
map as well as the parameter-to-prediction map around the maximum a
posterior point. We have found that for ice sheet
flow problems with the basal sliding parameter as the field of
interest, such linearizations are satisfactory
approximations for what would otherwise be an intractable problem
\cite{PetraMartinStadlerEtAl14}.

However, even with these linearizations, computing the covariance
of each of the resulting pdf's is prohibitive due to the need to solve
the forward ice sheet model a number of times equal to the parameter
dimension (or data dimension). We overcome this difficulty by
recognizing that these maps are inherently low-dimensional, since the
data inform a limited number of directions in parameter space, and the
predictions are influenced by a limited number of directions in
parameter space. Thus, with the right algorithm, the work---as
measured by ice sheet Stokes solves---should scale only with
the ``information dimension.'' The key idea to achieve this is to
construct a low rank approximation of the parameter-to-observable map
via a matrix-free randomized SVD method.

The result is a data-to-prediction framework whose computational cost
is overwhelmingly dominated by ice sheet model solves, both forward
and adjoint (and to a lesser extent elliptic solves representing
the action of parameter prior covariances). {\em Scalability of the
  entire data-to-prediction framework then follows when we show that
  (1) the number of forward (and adjoint) ice sheet solves needed for
  the data-to-prediction process is independent of the state dimension, the
  parameter dimension, and the data dimension; and (2) the forward
  (and adjoint) ice sheet solver demonstrates strong and weak scalability with
  increasing number of processor cores.}
We will show that our data-to-prediction framework---despite being
adamantly ``intrusive'' to ensure algorithmic scalability of the
inversion, uncertainty quantification, and prediction operations---can
be expressed in terms of a fixed and dimension-independent number of
forward-like ice sheet model solves, and thus exploits the same
algorithms, solver, and parallel implementation needed for the forward
problem. Thus, if a forward solver with both algorithmic and parallel
scalability can be designed---as will be shown in Section
\ref{sec:forward}---scalability of the entire data-to-prediction process
ensues.

We demonstrate this scalability of the data-to-prediction framework on
the problem of predicting the present-day ice mass flux from
Antarctica, starting
from Interferometric Synthetic Aperture Radar (InSAR) satellite
observations of the surface ice flow velocities, inferring the basal
sliding parameter field from this data via a 3D nonlinear
Stokes ice flow model on the present-day ice sheet geometry (Section
\ref{sec:inverse}), quantifying the 
uncertainty in the parameter inference using the Bayesian framework
(Section \ref{sec:bayesian}), and propagating the uncertain basal
sliding parameters through the forward ice sheet flow model to yield
predictions of ice mass flux with quantified uncertainties (Section
\ref{sec:uncertainty_propagation}).

\section{Forward problem: Modeling ice sheet flow}
\label{sec:forward}

The forward nonlinear Stokes ice sheet flow solver is the fundamental
kernel of our data-to-prediction framework and is invoked repeatedly
throughout. It is crucial that this forward solver scales
algorithmically and in parallel.  In this section we give a brief
summary of the design of our forward ice sheet flow solver and provide
performance results. For more details see~\cite{IsaacStadlerGhattas14}.

The flow of ice is commonly modeled as a viscous, shear-thinning,
incompressible fluid~\cite{Hutter83,Paterson94}. The balance of mass
and linear momentum state that
\begin{subequations}\label{eq:forward}
  \begin{align}
    -\gbf{\nabla} \cdot [\eta(\uu) (\gbf{\nabla u} + \gbf{\nabla
        u}^T) - \gbf I p] &= \rho
    \gbf{g},\label{eq:momentum}\\ \gbf{\nabla} \cdot \gbf{u} &=
    0, \label{eq:mass}
  \end{align}
where $\gbf u$ denotes the ice flow velocity, $p$ the pressure, $\rho$
the mass density of the ice, and $\gbf g$ the acceleration of gravity.
We employ a constitutive law for ice that relates the stress tensor
$\gbf \sigma$ and the strain rate tensor $\edot = \frac12 (\gbf{\nabla
  u} + \gbf{\nabla u}^T)$ by Glen's flow law~\cite{Glen55},
\begin{equation}
  \gbf \sigma = 2 \eta(\uu) \edot - \gbf I p, \; \text{with }
  \eta(\uu)
  =  \frac12 A^{-\frac1n} \; \secinve^{\frac{1-n}{2n}},\label{eq:glenslaw}
\end{equation}
where $\eta$ is the effective viscosity, $\gbf I$ is the second order
unit tensor, $\secinve = \frac12 \mathrm{tr}(\edot_{\gbf{u}}^2)$ is the
second invariant of the strain rate tensor, $n\ge 1$ is Glen's flow law
exponent, and $A$ is a flow rate factor that is a function of the
temperature $T$, parameterized by the Paterson-Budd
relation~\cite{PatersonBudd82}.  To construct a temperature field for
our simulations, we approximately solve steady-state, one-dimensional
advection-diffusion equations in the vertical direction at every horizontal
grid point. The advection velocity for this problem is only vertical, to
spatially decouple the columns.%
\footnote{%
  It should be noted that this temperature does not account for
  horizontal convection, and thus biases towards warmer temperatures at
  the ice sheet margins, where in reality colder ice from the sheet's
  interior is found at depth.%
} %
The advection velocity
for these equations interpolates between the accumulation
rate at the upper surface and zero at the base; we use the surface temperature
as the upper boundary condition and either the geothermal heat flux or the
temperature pressure melting point as the lower boundary condition where
appropriate.  We note that whether the resulting basal temperature is
below the pressure melting point does not affect our choice of boundary
conditions for the velocity described below; incorporating a regime
change between frozen and sliding basal conditions into our inversion
framework is the subject of future work.
The surface temperature, accumulation rate, and geothermal
heat flux data come from the ALBMAP dataset \cite{LeBrocqPayneVieli10}.

The top boundary of the Antarctic ice sheet $\Gamma_{\smallsub t}$ is
traction-free, and on its bottom boundary $\Gamma_{\smallsub b}$, we
impose no-normal flow and a Robin-type condition in the tangential direction:
\begin{align}
    \vec{\sigma} \vec{n} &= \vec{0} \quad \text{ on }
    \Gamma_{\smallsub t}, \label{eq:bcs1}\\
    \uu \cdot \vec n  = 0, \,\,\,\,
    \vec T \vec{\sigma}_{\!\uu} \vec{n} + \exp(\beta) \vec T\uu &= \vec
    0 \quad \text{ on } \Gamma_{\smallsub b}. \label{eq:bcs2}
\end{align}
\end{subequations}
Here, the coefficient $\exp(\beta)$ in the Robin boundary condition in
\eqref{eq:bcs2} describes the resistance to basal sliding. In the
following, we refer to $\beta = \beta(\bs x)$ as the basal sliding parameter
field.\footnote{Note that, in the literature, the parameterization of the Robin
  coefficient varies, e.g., instead of $\exp(\beta)$, sometimes $\beta^2$ is
  chosen in \eqref{eq:bcs2}.} Moreover, $ \gbf T :=
\gbf I - \gbf n \otimes \gbf n$ is a projection
operator onto the tangential plane, where ``$\otimes$'' denotes the
outer product. Note that $\exp(\beta)$, which relates tangential velocity to
tangential traction, subsumes several complex physical phenomena
and thus does not itself represent a physical parameter. It depends on
a combination of the frictional behavior of the ice sheet, the roughness of the
bedrock, the thickness of a plastically deforming layer of till,
and the amount and pressure of water present between the ice sheet and the
bedrock, all of which are presently poorly understood and constrained
sub-grid scale processes. As such, the basal sliding parameter field $\beta$
is subject to great uncertainty; subsequent 
sections will discuss the inverse problem of estimating it from
observational data. We note that in our model, the temperature does not
directly affect the boundary condition: a Robin condition is assumed
even when the temperature is below the pressure melting point.
On lateral ice-ocean boundaries, we use traction-free boundary
conditions above sea level. Below sea-level, the
normal component of the traction is set to the hydrostatic pressure
of sea water and the tangential components of the traction is zero.
We impose no-slip boundary conditions at
lateral boundaries, where the ice sheet terminates on land.
The lateral boundaries involve neither the observations nor the
uncertain basal parameters, and thus in the interests of conciseness,
we omit them from the discussion of the deterministic and statistical
inverse problems in Sections~\ref{sec:inverse} and \ref{sec:bayesian}.

The geometric description of the Antarctic ice sheet is constructed from the
ALBMAP dataset \cite{LeBrocqPayneVieli10}.  In our simulations, we
restrict ourselves to the
grounded portion of the ice sheet, i.e., we neglect ice
shelves, the extension of the sheet onto the surface of the
ocean: we therefore apply lateral boundary conditions at grounding lines
as discussed above.  A locally
refined mesh of hexahedral elements is used to discretize the ice sheet
domain.  We construct a coarse quadrilateral mesh that describes the lateral
geometry of the ice sheet and use this mesh as the basis for
forest-of-quadtree mesh refinement using the p4est library
\cite{BursteddeWilcoxGhattas11}.  Each quadrant in the refined mesh is then
used as the footprint of a column of hexahedra.  We do not constrain the
columns to have the same vertical resolution,
but allow each hexahedron to be
independently refined in the vertical direction.  This hybrid mesh refinement
strategy gives us flexibility to control the quality of the elements in our
computational mesh.  In particular, we can control the aspect ratio of the
elements using local refinement, which is not possible when isotropic
refinement is used, such as octree-based refinement. Note that as the
mesh is refined, we also refine the geometry description of the
Antarctic ice sheet.

For accuracy and efficiency,
we use a high-order accurate and locally volume-conserving
discretization that is provably stable for our locally and
nonconformingly\-/refined hexahedral meshes. To be precise, we use the
velocity/pressure finite element pair $\mathcal{Q}_k\times
\mathcal{Q}^{\mbox{\scriptsize disc}}_{k-2}$ 
for polynomial velocity order $k\ge 2$, i.e., with continuous
tensor-product polynomials of order $k$ for each velocity component,
and discontinuous tensor-product polynomials of order $k-2$ for the
pressure \cite{HeuvelineSchieweck07, ToselliSchwab03}.

An inexact Newton-Krylov method is used to solve the nonlinear systems
arising upon discretization of \eqref{eq:forward}, i.e., each Newton
linearization is solved inexactly using an iterative Krylov subspace
method. If the inexactness is properly controlled, the number of
overall Krylov iterations---and thus the overall work---is
minimized~\cite{IsaacStadlerGhattas14}.  Given a velocity/pressure
iterate $(\uu^j,p^j)$ at the $j$th Newton step, the new iterate is
computed as
$(\uu^{j+1},p^{j+1}) = (\uu^{j},p^{j}) + \alpha (\bs{\delta}_{\uu},\bs{\delta}_p)$
with an appropriate step length $\alpha>0$. The update
$(\bs{\delta}_{\uu},\bs{\delta}_p)^T$ is the solution of the linear
system 
\begin{equation}\label{eq:Newton}
  \begin{pmatrix}
    \bs{A}(\uu^j) & \bs{B}^T \\
    \bs{B} & \bs{0}
  \end{pmatrix}
  \begin{pmatrix} 
    \bs{\delta}_{\uu} \\ \bs{\delta}_p
  \end{pmatrix}=-
  \begin{pmatrix}
  \bs{r}_1^j \\
  \bs{r}_2^j
  \end{pmatrix},
\end{equation}
where $\bs{A}(\cdot)$ corresponds to the linearization of the
nonlinear viscous block in \eqref{eq:momentum} about $\uu^j$, $\bs{B}$
is the discretized divergence operator and $\bs{r}_1^j$ and
$\bs{r}_2^j$ are the nonlinear residuals corresponding to
\eqref{eq:momentum} and \eqref{eq:mass}, respectively, evaluated at
$(\uu^j,p^j)$. We use matrix-free matrix-vector applications to compute linear and
nonlinear residuals, which is particularly efficient for high-order
discretizations as it allows one to exploit the tensor product
structure of the finite element basis \cite{DevilleFischerMund02,
  BursteddeGhattasGurnisEtAl08, Brown10, IsaacStadlerGhattas14}.

Note that due to the form of the nonlinear rheology
\eqref{eq:glenslaw}, the linearization of the nonlinear viscous block
in \eqref{eq:momentum} results in a anisotropic tensor effective
viscosity in the operator $\bs{A}(\cdot)$. Due to self-adjointness,
this operator is identical to the viscous block operator of the
adjoint Stokes equations. The adjoint equations will be presented in
the next section; the expression for the anisotropic tensor effective
viscosity is given  in \eqref{eq:tensor-viscosity}.

To solve \eqref{eq:Newton} with an iterative Krylov method (GMRES or
FGMRES), efficient preconditioning is critical for scalability. We use
a block right preconditioner given by
\begin{equation*}
  \begin{pmatrix}
    \bs{A}(\uu^j) & \bs{B}^T \\
    \bs{0} & \tilde{\bs{S}}
  \end{pmatrix}^{-1},
\end{equation*}
where the inverse of the (1,1)-block $\bs{A}(\uu^j)$ is approximated
by an algebraic multigrid V-cycle (in particular, GAMG from PETSc
\cite{BalayBrownBuschelmanEtAl12}), and $\tilde{\bs{S}}$ denotes a
Schur complement approximation given by the inverse viscosity-weighted
lumped mass matrix \cite{ElmanSilvesterWathen05}.  Algebraic multigrid
for the (1,1)-block requires an assembled fine grid operator, but the
matrix $\bs{A}(\uu^j)$ is expensive to construct for high-order
discretizations: for our elements, the assembly time of
$\bs{A}(\uu^j)$ scales as $k^7$. The increasing density of
$\bs{A}(\uu^j)$ with $k$ also adds to the setup cost of the AMG
hierarchy, which is the bottleneck for parallel scalability.  We
therefore construct an approximation $\tilde{\bs{A}}(\uu^j)$ for
preconditioning based on low order finite elements \cite{Brown10}; in
the presence of nonconforming element faces, this can be challenging
\cite{IsaacStadlerGhattas14}.  The assembly time for
$\tilde{\bs{A}}(\uu^j)$ depends only on the total problem size, not on
the order $k$, and the sparsity is the same as for a low order
discretization.

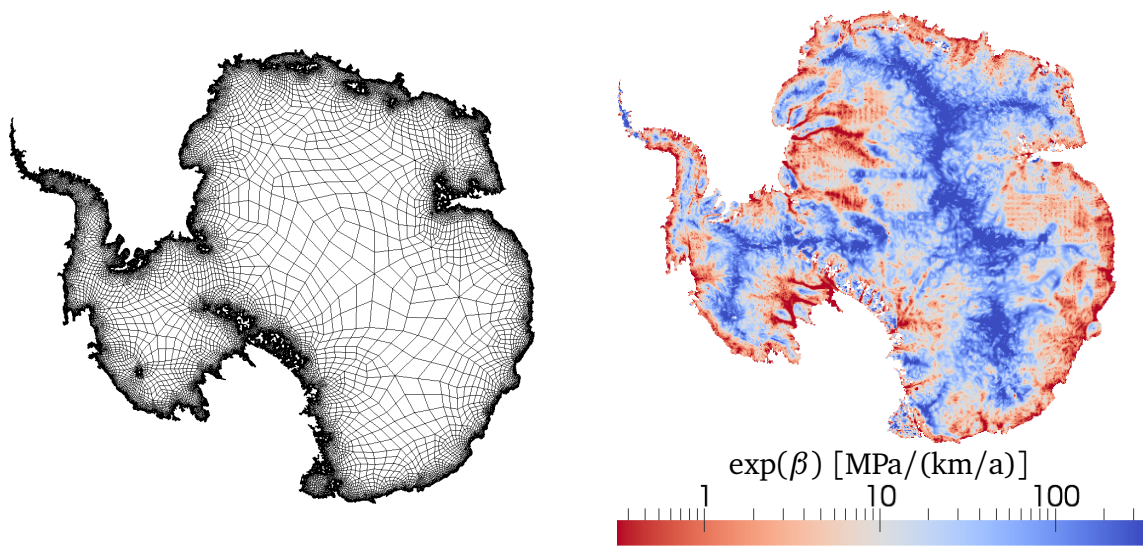
\begin{figure}[ht]\centering
\begin{tikzpicture}
  \node (l) at (-4,0.5)
    {\null\hbox{\includegraphics[width=.42\columnwidth]{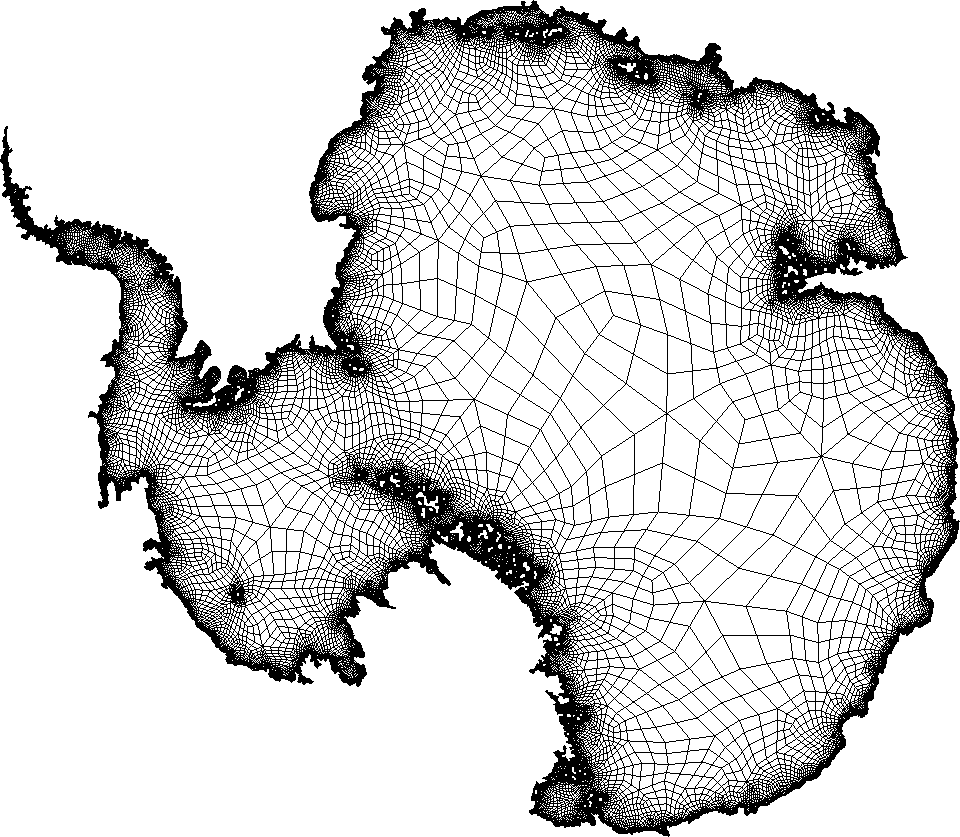}}};
  \node (r) at (4,0.3)
    {\null\hbox{\includegraphics[width=.42\columnwidth]{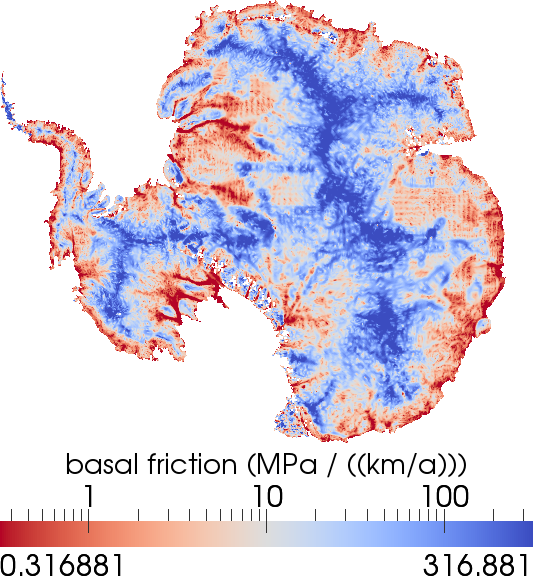}}};
    \draw[fill=white,white] (1,-2.2) rectangle (7,-1.8);
    \node at (4,-2) {\large$\exp(\beta)$ [MPa/(km/a)]};
    \draw[fill=white,white] (0,-3.7) rectangle (8,-3.1);
\end{tikzpicture}
  \caption{\label{fig:physics}
   Left: the coarse quadrilateral mesh from which the 3D
    meshes used in problems P1, P2, and P3 are extruded and refined.
    P1 is a discretization of the ice sheet with a maximum element
    width-to-height aspect ratio of 10:1; P2 and P3 are created by
    successive bisection in each direction. 
    Each problem uses $\mathcal{Q}_3 \times \mathcal{Q}^{\mbox{\tiny
        disc}}_1$ finite elements, with 330K, 2.6M, and 21M elements,
    respectively.  This amounts to 38M, 270M, and 2.1B unknowns.
    Right: the exponential of the synthetic basal sliding field $\beta$ used in
    the Robin-type basal boundary condition for problems P1, P2, and
    P3.
    The pre-exponential $A^{-1/n}$ in \eqref{eq:glenslaw} is computed
    from the Paterson-Budd relation \cite{PatersonBudd82}, with the
    temperature varying between $-50^{\circ}\mathrm{C}$ and
    $-10^{\circ}\mathrm{C}$.  For Glen's exponent, we use $n=3$, and
    we take $\rho = 910 kg/m^3$ and $g =9.81 m/s^2$.}
\end{figure}

Figure \ref{fig:physics} defines the Antarctic ice sheet problem we use
to study scalability of our nonlinear Stokes solver, described above.
The basal sliding parameter field was computed from relating the driving
stress due to gravity to the observed surface velocities; this synthetic
$\beta$ field is intended to be representative of the ``true'' $\beta$
field (one computed via inversion).
The base coarse mesh is shown, from which a sequence of
successively finer meshes is constructed. The figure also displays the
synthetic basal sliding parameter field $\beta$ used for the scaling study.
Table \ref{tab:titan} presents algorithmic and (strong) parallel scaling
studies on the sequence of successively refined Antarctic ice flow
problems defined in Figure \ref{fig:physics}. 
The results demonstrate
excellent algorithmic scalability, despite the severe computational
challenges of the problem (nonlinearity, indefiniteness, high-order
discretization, strongly-varying coefficients, anisotropy, high aspect
ratio mesh and geometry, locally-refined/nonconforming mesh, multigrid
for a vector system). Both the number of Newton iterations and the
average number of Krylov iterations per Newton iteration are only mildly
sensitive to problem size (over a range from 38M to 2.1B unknowns) and
number of cores (from 128 to 131,072 cores).
At larger core counts, the setup phase of algebraic multigrid remains
the bottleneck in parallel scalability.
Geometric multigrid for problems on
forest-of-octree meshes has proven efficient and scalable
\cite{SundarBirosBursteddeEtAl12}, and
can eliminate most of the setup time, but must be adapted to the
hybrid refinement scheme discussed above. 
This is a focus of our ongoing work.

Now that we have at our disposal a scalable forward nonlinear Stokes
solver, we proceed to the inverse problem and its
associated uncertainty, both of which require thousands of forward
Stokes solves, as will be seen in the next two sections. 

\begin{table}
  \begin{center}
    \renewcommand{\arraystretch}{1.1}
    \begin{tabular}{|@{\;}c@{\;}|@{\;}c@{\;}|@{\;}r@{\;}|@{\;}c@{\;}|@{\;}c@{\;}|@{\;}r@{\;}|@{\;}r@{\;}||@{\;}c@{\;}|}
      \hline
      & {\bf\#dof} & {\bf\#cores} & {\bf\#Newton} & \bf{\#Krylov} &
      \parbox[c][3em][c]{1.8cm}{\centering {\bf solve\\ time (s) / eff~(\%)}}
      &
      \parbox[c][3em][c]{1.8cm}{\centering {\bf setup\\ time~(s) / eff~(\%)}}
      &
      \parbox[c][3.5em][c]{1.4cm}{\centering {\bf \#Krylov\\ (Poisson)}}
      \\ \hline \hline

      \multirow{4}{*}{P1} & \multirow{4}{*}{38M } & 128    & 8  & 149 & 504.8 / 100 & 493.5  / 100     & 12 \\ \cline{3-8}
                          &                       & 256    & 8  & 153 & 259.6 / 97  & 260.4  / 95      & 12 \\ \cline{3-8}
                          &                       & 512    & 8  & 157 & 134.3 / 94  & 156.0  / 80      & 12 \\ \cline{3-8}
                          &                       & 1024   & 8  & 147 & 70.1  / 90  & 97.2   / 63      & 12 \\ \hline \hline
      \multirow{3}{*}{P2} & \multirow{3}{*}{270M} & 1024   & 9  & 240 & 796.6 / 100 & 735.0  / 100     & 12 \\ \cline{3-8} 
                          &                       & 2048   & 9  & 245 & 414.3 / 96  & 424.6  / 87      & 12 \\ \cline{3-8}
                          &                       & 8192   & 9  & 243 & 130.7 / 76  & 229.0  / 40      & 13 \\ \hline \hline
      \multirow{3}{*}{P3} & \multirow{3}{*}{2.1B} & 16,384  & 13 & 314 & 771.5 / 100 & 1424.5 / $\star$ & 15 \\ \cline{3-8}
                          &                       & 65,536  & 13 & 367 & 504.2 / 38  & 1697.1 / $\star$ & 15 \\ \cline{3-8}
                          &                       & 131,072 & 11 & 340 & 232.9 / 42  & 2033.1
      / $\star$ & 16 \\ \hline
    \end{tabular}
  \end{center}
  \caption{\label{tab:titan}
    Strong scaling results for our nonlinear Stokes solver for
    the Antarctic ice sheet problem on ORNL's Titan supercomputer. The
    table assesses algorithmic and parallel scalability for three
    problems (P1, P2, and P3) posed on successively finer
    meshes, each solved to a tolerance of $10^{-12}$.  A description
    of the problem setup is given in Figure \ref{fig:physics}. 
    For each run, we
    report the number of degrees of freedom ({\bf \#dof}), the CPU cores used ({\bf \#cores}), the number of
    Newton iterations ({\bf \#Newton}), the number of overall
    preconditioned Krylov iterations ({\bf \#Krylov}), and timings and
    parallel efficiencies for the matrix-vector multiplications and
    the multigrid V-cycles ({\bf solve time/eff}) and for the matrix
    assembly and the AMG setup ({\bf setup time/eff}).  A block-Jacobi/SOR
    smoother is used on each level of the V-cycle: incomplete
    factorization smoothers can be more effective for similar
    problems \cite{BrownSmithAhmadia13, IsaacStadlerGhattas14}, but
    it can be challenging to achieve robustness for the (1,1)-block
    for problems with high element aspect ratio, and strong variations
    in viscosity and bedrock topography. For
    comparison, we show the number of Krylov iterations ({\bf \#Krylov
    (Poisson)}) necessary to solve a scalar, constant coefficient,
    linear Poisson problem posed on the same meshes using
    $\mathcal{Q}_3$ finite elements and a V-cycle smoothed by damped
    block-Jacobi/IC(0), which yields close-to-optimal algorithmic
    scalability. The solution of this \emph{highly nonlinear,
      indefinite, highly heterogeneous, anisotropic, vector-valued}
    ice flow problem requires just
    $\sim 2 \times$ more Krylov iterations per Newton step than is
    needed for the (linear, positive definite, homogeneous, isotropic,
    scalar-valued) Poisson problem, indicating just how efficient our
    algorithms are. 
      The remaining bottleneck for scalability to very large
      core counts is the setup phase of the algebraic multigrid
      preconditioner. The ``$\star$'' symbol in the table indicates that the
      AMG setup times {\em increased} with strong scaling, and thus
      parallel efficiencies are not reported.}
\end{table}

\section{Inverse problem: Inferring the basal sliding parameter from observed
  ice surface velocity}
\label{sec:inverse}

In this section, we describe the solution of the inverse problem of
inferring the two-dimensional basal sliding parameter field $\beta$ in
\eqref{eq:bcs2} from satellite-derived surface velocity
observations with the ice sheet flow model described in the previous
section. Several recent studies have focused on inversion for the
basal boundary conditions in ice sheet models.  Contributions that use
adjoint-based gradient information to solve the deterministic inverse
problem include~\cite{VieliPayne03, JoughinMacAyealTulaczyk04,
  LarourRignotJoughinEtAl05, MorlighemRignotSeroussiEtAl10,
  GoldbergSergienko11,
  LarourSeroussiMorlighem12, PeregoPriceStadler14,
  GoldbergHeimbach13}. In~\cite{MorlighemSeroussiLarourEtAl13}
the authors use InSAR surface velocity measurements, as employed in
this section, and invert for the basal sliding parameter field for the
Antarctic ice sheet. However, unlike the present work, the authors use
the first-order accurate Stokes approximation~\cite{Blatter95, Pattyn03}. This model is derived from the full
Stokes equations by making the assumptions that horizontal gradients
of vertical velocities are negligible compared to vertical
gradients of horizontal velocities, that 
horizontal gradients of vertical shear stresses are 
small compared to $\rho g$, and that the pressure in
vertical direction is hydrostatic.
The work presented here represents the first continental scale ice
sheet inversion using the full Stokes ice flow equations, which is
considered the highest fidelity ice flow model available. Moreover,
this work is the first to quantify uncertainty in the solution of the
Antarctic ice sheet inverse problem using any flow model (Section
\ref{sec:bayesian}), as well as propagate that uncertainty to a
prediction quantity of interest (Section
\ref{sec:uncertainty_propagation}).

In the remainder of this section, we present the inverse problem
formulation, give a brief overview of the solution method, provide
expressions for the gradient based on adjoint ice flow equations, and
present results for the deterministic inverse problem.

\subsection{The regularized inverse problem}
\label{detip}
The inverse problem is formulated as follows: given (possibly noisy)
observational data $\uobs$ of the ice surface velocity field,
we wish to infer the basal sliding parameter field $\beta(\bs{x})$
(defined as a coefficient in a Robin boundary condition at the base of
the ice sheet; see \eqref{eq:bcs2}) that produces a surface velocity field that best fits
the observed velocity.  This can be formulated as the variational
nonlinear least squares optimization problem
\begin{align}
\label{objfunction} 
  \min_{\beta} \mathcal{J}(\beta) := \; \frac12
  \int_{\Gamma_{\smallsub t}} \frac{|\mathcal{B}\uu(
  \beta)-\uobs|^2}{|\uobs|^2 + \varepsilon} \; d\gbf{s} +
  \Reg(\beta),
\end{align}
where the forward velocity $\uu$ is the solution of the forward
nonlinear Stokes problem~\eqref{eq:momentum}--\eqref{eq:bcs2} for a
given basal sliding parameter field $\beta$, and $\B$ is an observation
operator that restricts the model ice velocity field to the top surface.
Since the observed surface flow velocities $\uobs$ can vary over three
or more orders of magnitude, we normalize the data misfit term by
$|\uobs|^2 + \varepsilon$, where $\varepsilon$ is a small constant
to prevent division by zero. 

The regularization term 
\begin{equation}\label{eq:reg-form}
\mathcal R(\beta) := \frac12
\|\mc{A}^\kappa(\beta-\beta_0)\|_{L^2(\Omega)}^2
\end{equation}
penalizes oscillatory
components of the basal sliding parameter field on the basal surface of the
ice, thus restricting the solution to smoothly varying 
fields. Here $\beta_0$ is a reference basal sliding parameter, the
differential operator $\mc{A}$ is defined as
\begin{equation}\label{eq:reg}
  \mc{A}(\beta):=
  \left\{
  \begin{array}{rrl}
    - \gamma \gbf \Delta_{\Gamma} \beta + \delta \beta
    \quad &&\text{in } \Gamma_{\smallsub b},\\
    (\gamma \gbf \nabla_{\Gamma} \beta) \cdot \gbf \nu \quad &&\text{on }
    \partial \Gamma_{\smallsub b},
  \end{array}\right.
\end{equation}
where $\gamma > 0$ is the regularization parameter that controls the
strength of the imposed smoothness relative to the data misfit,
$\delta > 0$, typically small compared to $\gamma$, is added to make the
regularization operator invertible, $\gbf \nabla_{\Gamma}$ is the
tangential gradient, $\gbf \Delta_{\Gamma}$ is the Laplace-Beltrami
operator,
and $\gbf \nu$ denotes the outward unit normal on $\partial \Gamma_{\!b}$.
The regularization operator is thus a positive definite elliptic
operator of order $4\kappa$. The need for such a regularization stems from
the fact that small wavelength components of the basal sliding parameter field
cannot be inferred from surface velocity observations due to the
smoothing nature of the map from $\beta$ to $\mathcal{B}\uu( \beta)$
\cite{PetraZhuStadlerEtAl12}. In the absence of such a term, the
inverse problem is ill-posed, that is, its solution is not unique and
is highly sensitive to errors in the observations (for general
references on regularization in inverse problems, see e.g.,
\cite{EnglHankeNeubauer96,Vogel02}).

When solving the deterministic inverse problem, a value of $\kappa =
\frac12$ within the regularization operator \eqref{eq:reg-form} (which
corresponds to classical $H^1$-type Tikhonov regularization) suffices
for well-posedness. On the other hand, in the Bayesian inference
framework (discussed in Section~\ref{sec:bayesian}), in which the
regularization term \eqref{eq:reg-form} reflects prior knowledge on the
basal sliding parameter field, higher orders of the elliptic
regularization operator are required for well-posedness, depending on
the spatial dimension of the domain $\Omega$~\cite{Stuart10}. For
example, sufficient values are $\kappa = \frac12$ in one dimension and
$\kappa = 1$ in two and three dimensions.  In this setting,
$\mc{A}^{-\kappa}$ is the covariance operator of the prior distribution.
The Green's function of $\mc{A}^{\kappa}$ at a point $x$ describes
the correlation between the parameter value at $x$ and values elsewhere.
Our choice of differential operator makes this correlation decay
smoothly with increasing distance.  In practice, this allows us to apply
a dense spatial correlation operator by inverting $\mc{A}^{\kappa}$
without storing a dense matrix.

\subsection{Inverse problem solver: Inexact Newton-CG}
\label{ipmethod}

We solve the optimization problem \eqref{objfunction} with an inexact
matrix-free Newton--CG method. This amounts to (approximately) solving
the linear system that arises at each Newton iteration,
\begin{equation}\label{eq:newton-opt}
\mathcal{H}(\beta^k) \, \Delta \beta = - \mathcal{G}(\beta^k),
\end{equation}
by the conjugate gradient method, and then updating $\beta$ by 
\begin{displaymath}
\beta^{k+1} := \Delta \beta + \alpha^k \beta^k.
\end{displaymath}
Here, $\mathcal{G}(\beta)$ and $\mathcal{H}(\beta)$ are the gradient
and Hessian, respectively, of $\mathcal{J}(\beta)$, both with respect
to $\beta$, and $\alpha^k$ is a step length chosen by a suitable line
search method. The next section presents an efficient method for
computation of the gradient via a variational adjoint method; in
particular, the gradient is given by expression \eqref{eq:gradient},
which depends on solutions of forward and adjoint Stokes
problems. Hessian computations are made tractable by recalling that
the CG solver does not require the Hessian operator by itself; it
requires only the action of the Hessian in a given direction.
Operator-free (i.e., matrix-free after discretization) computation of
this Hessian action (i.e., Hessian-vector product) is presented in
Section \ref{hess}, in particular in the expression
\eqref{eq:hessian_action}; it involves second order adjoints but has a
structure similar to the gradient computation.

Details of the inexact Newton-CG method are provided in
\cite{PetraZhuStadlerEtAl12} and references therein. A summary of the
components of the method is as follows:
\begin{itemize}
\zapspace
\item The Newton system is solved inexactly by early termination of CG
  iterations via Eisenstat--Walker (to prevent oversolving) and
  Steihaug (to avoid negative curvature) criteria.

\item Preconditioning is effected by the inverse of the (elliptic)
  regularization operator, which is carried out by multigrid
  V-cycles.

\item Globalization is by an Armijo backtracking line search. 

\item Continuation on the regularization is carried out to
  warm-start the Newton iterations, i.e., we initially use a large
  value of $\gamma$ in \eqref{eq:reg} and decrease $\gamma$ during the
  iteration to the desired value.

\item As elaborated in Sections \ref{grad} and \ref{hess}, gradients
  and Hessian actions at each CG iteration are expressed in terms of
  solutions of forward and adjoint PDEs, and their linearizations.

\item Parallel implementation of all components, whose cost is dominated
  by solution of forward and adjoint PDEs and evaluation of inner
  product-like quantities to construct gradient and Hessian action
  quantities~\cite{AkcelikBirosGhattasEtAl06a}.

\end{itemize}
\zapspace

Parallel and algorithmic scalability follow as a result of the
following properties. Because the dominant components of the method
can be expressed as solutions or evaluations of Stokes-like systems
and inner products, parallel scalability---that is, maintaining high
parallel efficiency as the number of cores increases---is assured
whenever a scalable solver for the underlying PDEs is available (as
demonstrated in Section \ref{sec:forward}).
The remaining ingredient needed to obtain overall scalability is that
the method exhibit algorithmic scalability---that is, the number of
iterations does not increase with increasing problem size.  This is
indeed the case: for a wide class of nonlinear inverse problems,
the number of outer Newton iterations and of inner CG iterations is
independent of the mesh size and hence parameter dimension, as will be
demonstrated in Section \ref{ipres}.
This is a consequence of the use of a Newton solver, of the compactness of the
Hessian of the data misfit term in~\eqref{objfunction}, and the
preconditioning by the inverse of the regularization operator so that the
resulting preconditioned Hessian is a compact perturbation of the identity,
for which Krylov subspace methods exhibit mesh-independent iterations
\cite{CampbellIpsenKelleyEtAl96}.

\subsection{Gradient computation via a variational adjoint method}
\label{grad}

The Newton-CG method described in the previous section requires
computation of the (infinite dimensional) gradient
$\mathcal{G}(\beta)$, which is the Fr\'{e}chet derivative of
$\mathcal{J}(\beta)$ with respect to $\beta$. Using the Lagrangian
formulation and variational calculus, we can derive an expression for
$\mathcal{G}(\beta)$ at an arbitrary point $\beta$ in parameter space
\cite{PetraZhuStadlerEtAl12}; its strong form is given by
\begin{equation}
\label{eq:gradient}
  \mc{G}(\beta):= \exp(\beta)\; \gbf T \uu \cdot
  \gbf T \vv + \mc{A}^{2\kappa}(\beta-\beta_0) \quad \text{ on
  } \Gamma_{\smallsub b},
\end{equation}
where the forward velocity/pressure pair ($\uu$, $p$) satisfies the
forward Stokes problem~\eqref{eq:forward}, and the adjoint
velocity/pressure pair ($\vv$, $q$) satisfies the {\it adjoint Stokes
  problem},
\begin{subequations}\label{eq:adjoint}
\begin{alignat}{2}
  \gbf{\nabla} \cdot \vv &= 0  &\; &  \quad \text{ in }  \Omega,\\
  - \gbf{\nabla} \cdot \gbf{\sigma}_{\!\vv} &= \gbf{0} & & \quad
   \text{ in } \Omega,\\ 
  \gbf{\sigma}_{\!\vv} \gbf{n} &= - \mathcal{B}^*\bs\Gamma(\mathcal{B}\uu -\uobs
  &)\; & \quad \text{ on }
  \Gamma_{\!\mbox{\tiny} t}, \\
  T\gbf{\sigma}_{\!\vv} \gbf{n} + \exp(\beta) \gbf T\vv &=
  \gbf{0}, \;\; \vv\cdot \gbf{n} = 0 &\; & \quad \text{ on }
  \Gamma_{\!\mbox{\tiny} b}. 
\end{alignat}
\end{subequations}
Here, $\bs\Gamma$ accounts for the pointwise scaling by
$(|\uobs|^2+\varepsilon)^{-1}$ in \eqref{objfunction}, and
 the adjoint stress $\gbf{\sigma}_{\!\vv}$ is given in terms of
the adjoint strain rate $\edot_{\vv}$ by
\begin{equation}\label{eq:tensor-viscosity}
  \gbf{\sigma}_{\!\vv} := 2\eta(\uu)
  \, \bigl(\mathsf{I} + \frac{1-{n}}{{n}} \frac{\edot_{\uu}\otimes
    \edot_{\uu}}{\edot_{\uu} \cdot \,
    \edot_{\uu}}\bigr)\edot_{\vv} -\gbf{I} q,
\end{equation}
where $\mathsf{I}$ is the fourth order identity tensor. 
The adjoint Stokes problem has several notable properties: (1) its
only source term is the misfit between observed and predicted surface
velocity on the top boundary; (2) while the forward problem is a
nonlinear Stokes problem, the adjoint Stokes problem is linear in the
adjoint velocity and pressure, and is characterized by a linearized
Stokes operator with a 4th-order tensor anisotropic effective
viscosity with anisotropic component that acts in the direction of the
forward strain rate $\edot_{\uu}$;
and (3) the adjoint Stokes operator is the same operator as the
Jacobian that arises in Newton's method for solving the forward
(nonlinear) Stokes problem (as a consequence of its
self-adjointness). Because of this equivalence, we use the same
discretization method and mesh, linear solver, and preconditioner for
the adjoint problem as was presented for (a Newton step of) the
forward problem in Section \ref{sec:forward}. Use of the same
discretization guarantees so-called {\em gradient consistency}
\cite{PetraZhuStadlerEtAl12}. 

In summary, to compute the gradient \eqref{eq:gradient} for a given
$\beta$ field, we first solve the forward nonlinear Stokes problem
\eqref{eq:forward} for the forward velocity/pressure pair ($\uu$,
$p$), and then using this forward pair, we solve the adjoint problem
\eqref{eq:adjoint} for the adjoint velocity/pressure pair ($\vv$,
$q$). Finally, both pairs of fields are used to evaluate the gradient
expression \eqref{eq:gradient} at any point on the basal surface
$\Gamma_{\smallsub b}$ or its boundary $\partial
\Gamma_{\smallsub b}$. The action of the Hessian in a given
direction could be computed by directional differences of gradients,
using the method just described to compute the gradient (necessitating
a nonlinear forward Stokes solve for each gradient
evaluation). However, it is more efficient to derive expressions for
the Hessian action in infinite-dimensional form and compute with them
since only a linearized forward solve is required for each Hessian
action; presentation of these expressions will be deferred to Section
\ref{hess}, where the Hessian plays an important role in
characterizing the uncertainty in the solution of the inverse problem.

\subsection{Inversion results}
\label{ipres}
Here we assess the performance of the inversion algorithm of the
previous section, and present continental-scale inversion results for
Antarctica.
The observational data are the Antarctic ice sheet surface ice
velocities obtained using satellite synthetic aperture radar
interferometry (InSAR)~\cite{RignotMouginotScheuchl11}. We initialize
the Newton iteration with a constant basal sliding parameter
field, whose value is such that the ice is strongly connected
to the bedrock and thus little sliding occurs.
The iteration is terminated as soon as the norm of the gradient of
$\mathcal{J}$ is decreased by a factor of $10^5$.
First we study the performance of the inexact Newton-CG method
described in the previous section as the parameter and data dimensions
that characterize the inverse problem grow. Both of these dimensions
are tied to the mesh size, since both the basal sliding parameter and
observational surface velocity are treated as continuous fields. To
make this study tractable, the inversion domain is limited to the Pine
Island Glacier region. Performance is characterized in terms of the
number of linear or linearized Stokes systems that must be solved,
since these are the core kernels underlying objective function,
gradient, and Hessian-vector product computations, and thus
overwhelmingly dominate the run time (the cost of the remaining linear
algebra associated with the Newton-CG optimization method is
negligible relative to that of the Stokes solves).

\begin{figure}[ht]
\begin{center}
\includegraphics[width=.4\columnwidth]{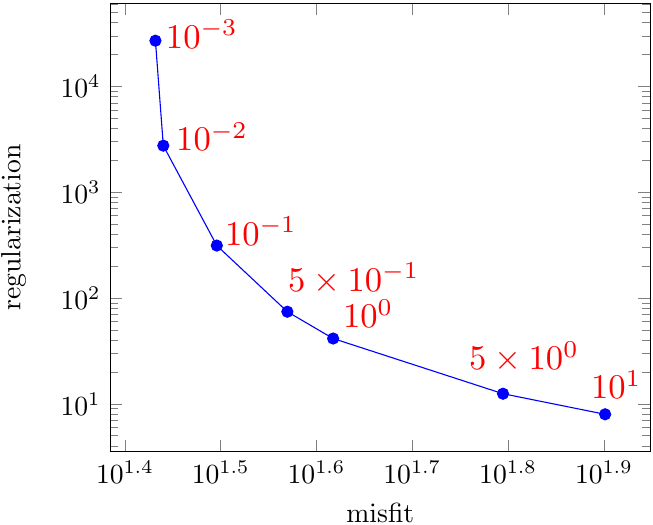}
\end{center}
\caption{L-curve-based regularization parameter selection. Each point
  (in blue) on the curve represents the solution of an inverse problem
  for a different choice of the value of the regularization parameter
  $\gamma$ (shown in red). For each value of $\gamma$, the vertical
  axis plots the magnitude of the regularization term
  $\mathcal{R}(\beta)$
  evaluated at the optimum value of $\beta$; the horizontal axis plots
  the value of the data misfit term i.e., the first term in
  $\mathcal{J}(\beta)$ defined by \eqref{objfunction}, also evaluated
  at the optimum value of $\beta$. The L-curve criterion takes as the
  ``optimum'' regularization parameter the value at the point of
  maximum curvature, in this case $\gamma \approx 0.1$. The criterion
  posits this value of $\gamma$ as the best tradeoff between
  minimizing the data misfit and controlling small-scale
  variations due to data noise in the
  resulting $\beta$ field.}
\label{fig:lcurve}
\end{figure}

Table~\ref{tbl:performance} presents algorithmic performance for a
sequence of increasingly finer meshes, and hence increasing state,
parameter, and data dimensions. The
regularization operator for these cases is the Laplacian, i.e., in the
operator $\mc{A}$ defined by \eqref{eq:reg}, we take $k = \frac12$,
with $\delta = 0$. We also set the value of the reference basal
sliding parameter field $\beta_0$ in \eqref{eq:reg-form} to 0.
An L-curve criterion is employed to determine an optimal value
of the regularization parameter. The L-curve criterion requires
solution of several inverse problems with different values of
$\gamma$.
Figure~\ref{fig:lcurve} depicts the L-curve for the inverse problem
defined by the first row of Table~\ref{tbl:performance}. Based on this
criterion, the regularization parameter for all tests is taken to be
$10^{-1}$. As can be seen in Table~\ref{tbl:performance}, the number
of Newton and CG iterations required by the inexact Newton-CG inverse
solver is insensitive to the parameter and data dimensions, when
scaling from 10K to 1.5M inversion parameters (and 96K to 23M states).
Thus, the number of Stokes solves does not increase with increasing
inverse problem size, leading to a perfectly scalable inverse solver.

\begin{table}
\begin{center}
  \small
  \begin{tabular}{|r|r|r|r|r|r|}
    \hline
        {\bf \#sdof}    &   {\bf \#pdof}   &    {\bf  \#N}   &  {\bf  \#CG}
        & {\bf  avgCG}   &  {\bf \#Stokes} \\
        \hline
        95,796 &  10,371 &  42 &  2718 &    65         &     7031 \\
        233,834 &  25,295 &  39 &  2342 &    60        &     6440 \\
        848,850 &  91,787 &  39 &  2577 &    66        &     6856 \\
        3,372,707 & 364,649 &  39 &  2211 &    57       &     6193 \\
        22,570,303 & 1,456,225 &  40 &  1923 &    48     &     5376\\
        \hline
  \end{tabular}
\end{center}
\vspace{0.05in} 
\caption{Algorithmic scalability of inverse solver for the ice sheet
  inverse problem posed on the Pine Island Glacier region of
  Antarctica. Number of Newton
  and CG iterations and number of linearized Stokes solves for a
  sequence of increasingly finer meshes (and hence inversion
  parameters) for the inexact Newton-CG method. 
  The first column ({\bf \#sdof}) shows the number of
  degrees of freedom for the state variables; the second column ({\bf
    \#pdof}) shows the number of degrees of freedom for the inversion
  parameter field; the third column ({\bf \#N}) reports the number of
  Newton iterations; columns four ({\bf \#CG}) and five ({\bf avgCG})
  report the total and the average (per Newton iteration) number of CG
  iterations;
and the last column ({\bf \#Stokes}) reports the total number of
linear(ized) Stokes solves (from forward, adjoint, and incremental
forward and adjoint problems). The Newton iterations are terminated
when the norm of the gradient is decreased by a factor of $10^5$ .
The CG iterations are terminated when the norm of the residual of the
Newton system drops below a tolerance that is proportional to the norm
of the gradient, per the Eisenstat--Walker criterion corresponding to ``Choice
2'' in \cite[Section 2]{EisenstatWalker96} (with safeguards as described in
Section 2.1 of the same).
These results illustrate that the cost of solving the inverse problem
by the inexact Newton-CG method, measured by the number of Stokes
solves, is independent of the parameter dimension. The cost is also
independent of the data dimension, since the surface observational
velocity field is refined with decreasing mesh size, just as the basal
sliding parameter field is refined. \label{tbl:performance}} 
\end{table}

\begin{figure}[ht]
\centering
\begin{tikzpicture}
  \node (l) at (0,0.5)
  {\includegraphics[width=.6\columnwidth]{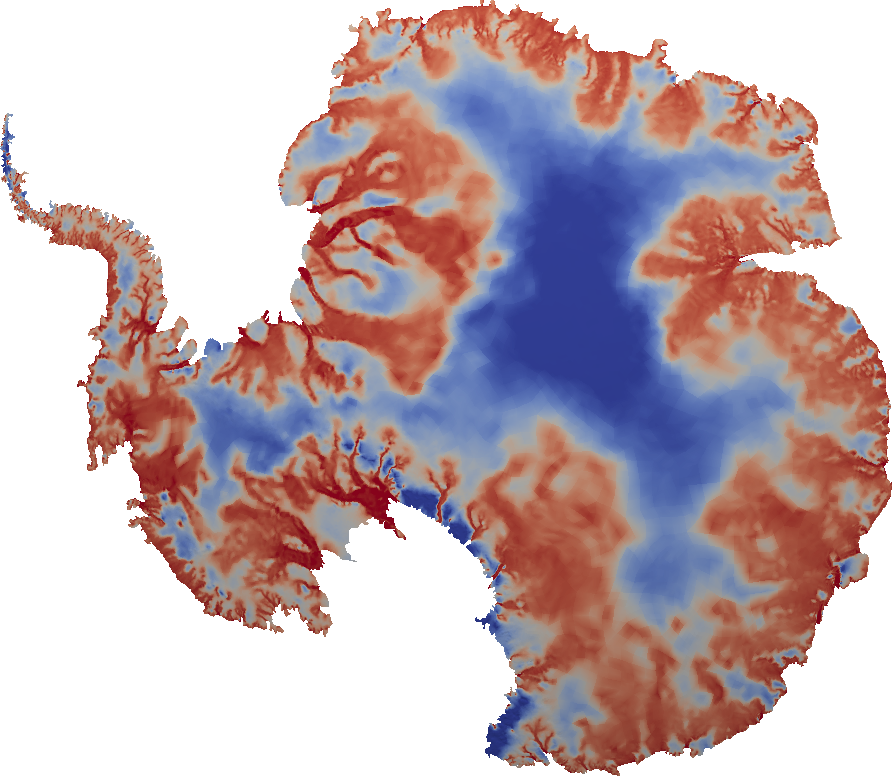}};
  \node (r) at (6,0.5)
  {\includegraphics[width=0.13\columnwidth]{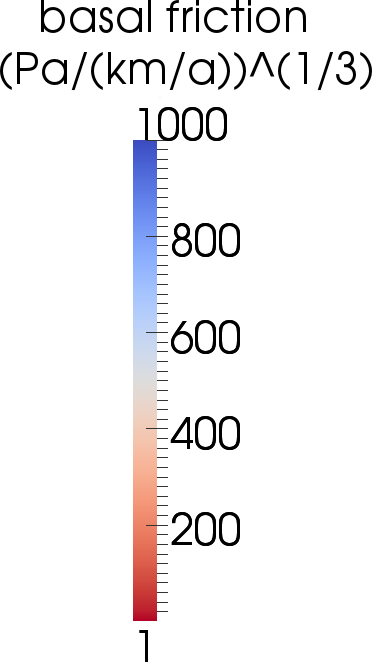}};
  \draw[fill=white,white] (4.8,1.8) rectangle (7.2,2.5);
  \node at (6,2.5) {$\exp(\beta)^{1/3}$};
  \node at (6,2.05) {[(Pa/(km/a))$^{1/3}$]};
\end{tikzpicture}
\caption{Solution of inverse problem of inferring the Antarctic basal
  sliding parameter field from InSAR surface velocity
  observations. For visualization purposes, the cube root of the
  exponential of the basal sliding coefficient is plotted; the actual
  Robin coefficient $\exp(\beta)$ of the basal sliding boundary
  condition varies over nine orders of magnitude. Low (red) and high
  (blue) values for the basal sliding parameter represent low and high
  resistance to sliding, and correlate with fast and slow ice flow
  regions, respectively.
}
\label{fig:map}
\end{figure}
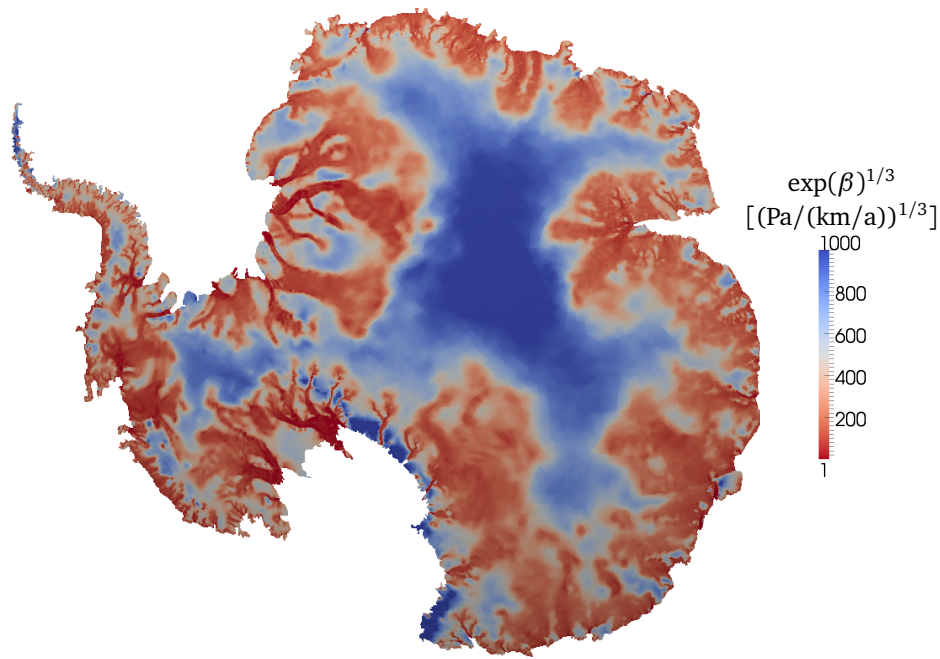

Figure~\ref{fig:map} depicts the solution of the ice sheet inverse
problem, i.e., the reconstruction of the Antarctic basal sliding parameter
field $\beta$ from InSAR observations of surface ice flow velocity and
the nonlinear Stokes ice flow model. Note that the basal sliding parameter field
varies over nine orders of magnitude.  Low (red) and high (blue)
values of  $\beta$ represent low  and high resistance to basal
sliding and correlate with fast and slow ice flow, respectively. As
can be seen, weak resistance to basal sliding extends deep into the
interior of the continent.

\begin{figure}[ht]\centering
\begin{tikzpicture}
\node (u) at (3,5)
{\includegraphics[width=.7\columnwidth]{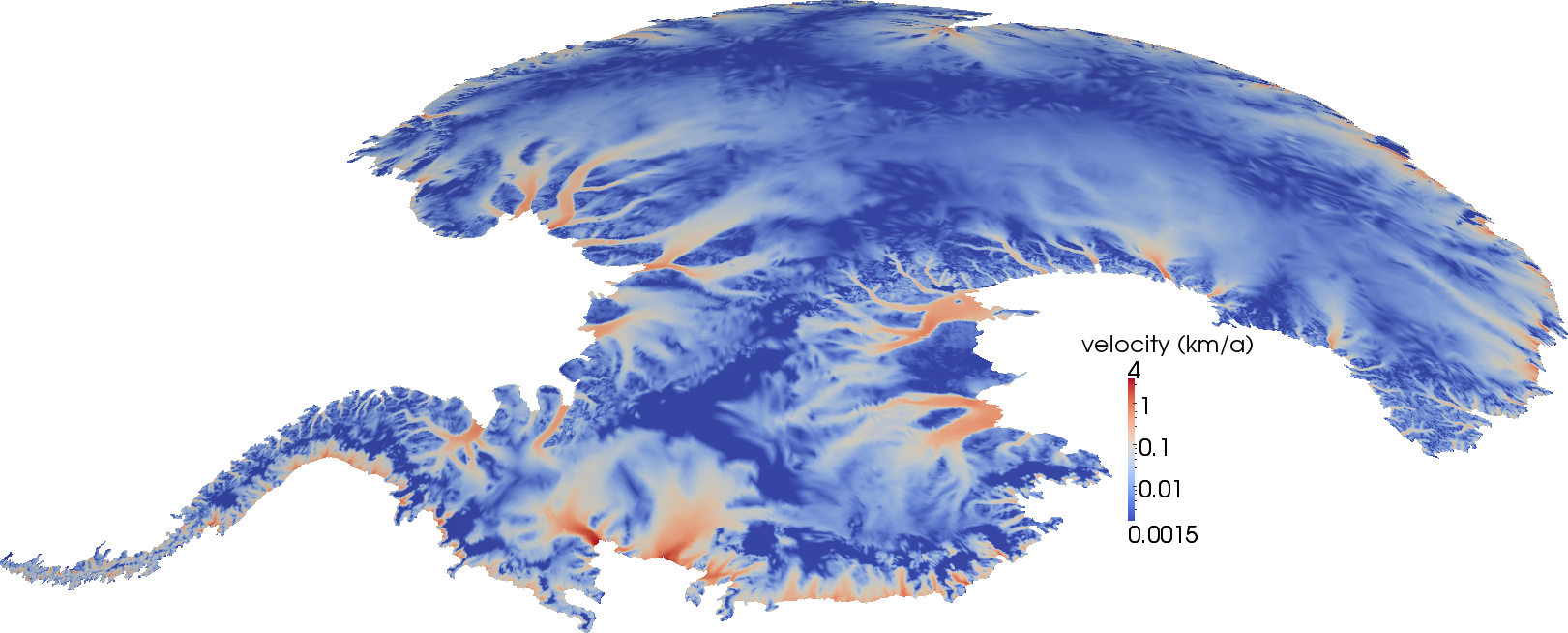}};
\node (u) at (3,0)
{\includegraphics[width=.7\columnwidth]{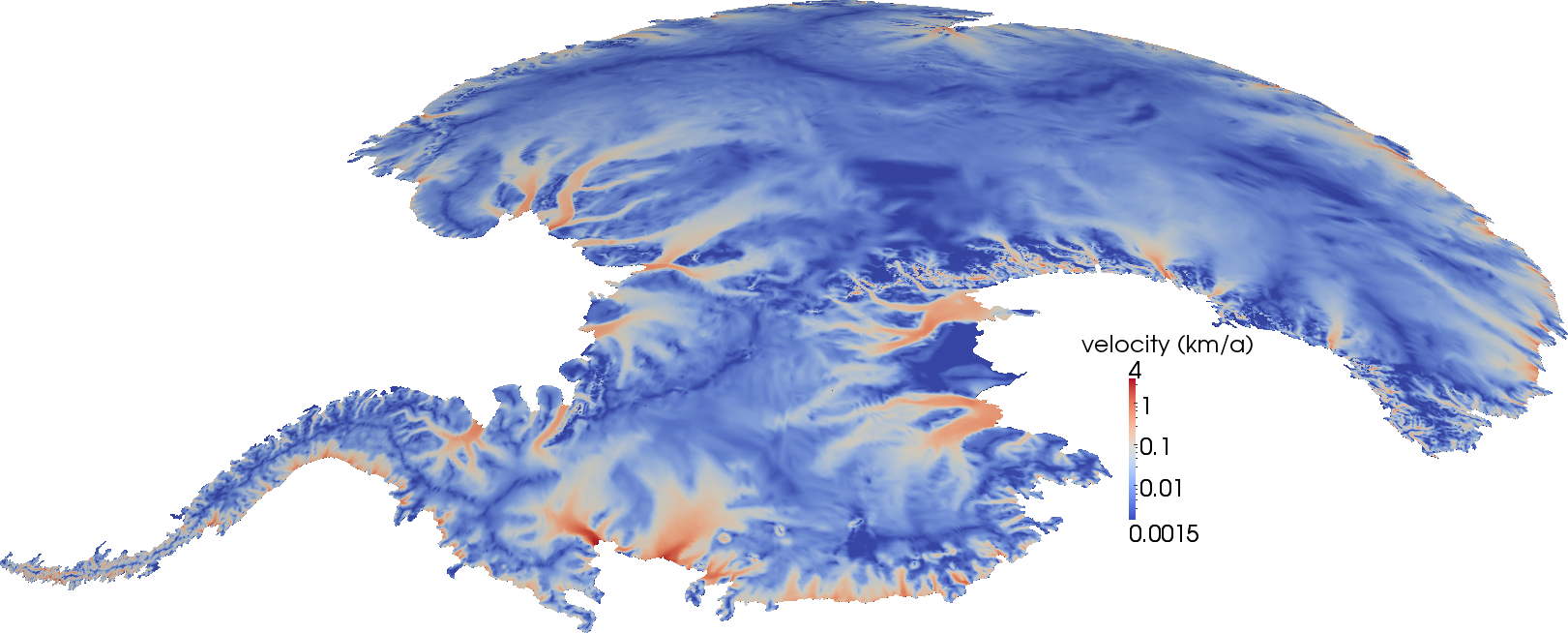}};
\node at (5.3,5.1){\textcolor{red!60!black}{\small Ross}};
\node at (0.2,5){\textcolor{red!60!black}{\small Filchner-Ronne}};
\node at (4.7,7.4){\textcolor{red!60!black}{\small Amery}};
\node at (2.2,2.8){\textcolor{red!60!black}{\small George VI}};
\node at (7.8,6.6){\textcolor{red!60!black}{\small Shackleton}};
\end{tikzpicture}
\caption{Observed (top) and reconstructed via solution of the inverse
  problem (bottom) surface velocity fields. In the top image, we also
  highlight the names and locations of the largest ice shelves, the
  extensions of the ice sheet onto the ocean. Most of the mass
  from the fast flowing ice streams contributes to ice shelves,
  and this mass is, over time, released to the ocean due to melting and
  iceberg calving.}
\label{fig:velocity}
\end{figure}
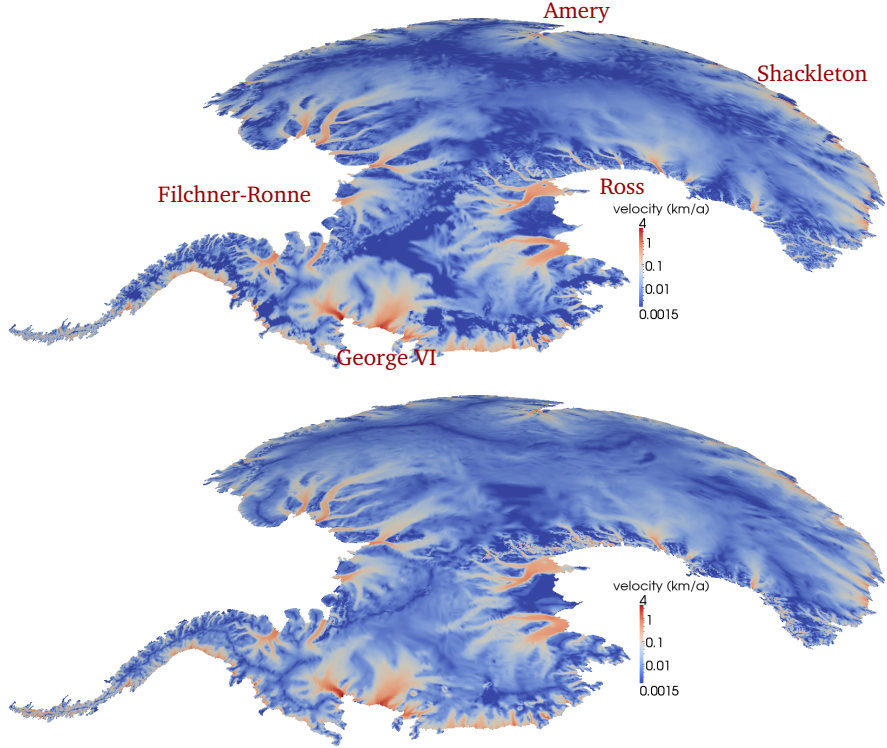

Figure~\ref{fig:velocity} addresses the question of how
successful the inversion is in creating a model that explains the
data. The top image portrays the observed surface velocity field over
the continent, which varies widely from a few meters per year in the
interior of the continent (in dark blue), to several kilometers per
year in the fast flowing ice streams and outlet glaciers (in bright
red). The bottom image shows the surface velocity field that has been
reconstructed from solution of the nonlinear Stokes ice flow model
using the inferred basal sliding parameter field. By visual inspection, we can see
that the two surface velocity fields agree well, particularly in fast
flow regions.
The difference between the two images reflects the ill-posed nature of
the inverse problem, for which data noise and model inadequacy can
amplify errors in poorly-inferred features of the inversion parameter
field. Nevertheless, the inversion is successful in inferring a basal
sliding parameter field and resulting ice sheet model that is able to
fit the data well, particularly in the critical regions of the ice sheet
that impact sea level (i.e., the fast flowing ice streams of marine ice
sheets that deliver the bulk of the mass flux to the ocean, and are most
sensitive to forcing changes at ice-ocean interfaces
\cite{JoughinAlleyHolland12}).

Despite this success in fitting the data, ultimately we are still
dealing with an ill-posed inverse problem in which the ``true
solution'' cannot be found exactly. The question we must address is:
what confidence do we have in the inverse solution we have obtained?
The deterministic inverse problem described in this section is not
equipped to answer this question. Instead, we turn to the framework of
Bayesian inference, which provides a systematic means of quantifying
uncertainty in the solution of the inverse problem. In the next
section we present algorithms for making the Bayesian framework
tractable for large-scale ice sheet inverse problems.

\section{Bayesian quantification of parameter uncertainty:
  estimating the posterior pdf of the basal sliding parameter field} 
\label{sec:bayesian}

In this section we tackle the problem of quantifying the uncertainty
in the solution of the ice sheet inverse problem discussed in the
previous section. We adopt the framework of Bayesian inference
\cite{Tarantola05, KaipioSomersalo05}, and in particular its extension
to infinite-dimensional inverse problems~\cite{Stuart10}. To keep our
discussion compact, we begin by presenting finite-dimensional
expressions (i.e., after discretization of the parameter space) for
the Bayesian formulation of the inverse problem; we refer the reader
to \cite{PetraMartinStadlerEtAl14} for elaboration of the
infinite-dimensional framework and associated discretization issues
for the ice sheet inverse problem. Next we discuss a low-rank-based
approximation of the posterior covariance (built on a low-rank
approximation of the Hessian of the data misfit) that permits %
very large parameter dimensions. We then present expressions that
show how the low-rank approximation of this Hessian can be computed
for the ice sheet inverse problem, and discuss properties of the
Hessian that, in combination with the proposed algorithm, permit
computation of the uncertainty in the inverse solution in a fixed
number of forward/adjoint solves, independent of the parameter
dimension. Finally, the methodology is applied to the large-scale
Antarctic ice sheet inverse problem.

\subsection{Bayesian solution of the inverse problem}
In the Bayesian formulation, we state the inverse problem as a problem of
\emph{statistical inference} over the space of uncertain parameters,
which are to be inferred from  data and a physical model.
The resulting solution to the statistical inverse problem is a
posterior distribution that assigns to any candidate set of parameter
fields our belief (expressed as a probability) that a member of this
candidate set is the ``true'' parameter field that gave rise to the
observed data.
When discretized, this infinite-dimensional inverse problem leads to a
large-scale problem of inference over the discrete parameters
$\vec{\beta} \in \mathbb R^n$, corresponding to degrees of freedom in
the parameter field mesh. 

The posterior probability distribution combines the prior pdf
$\pi_{\text{prior}}(\vec{\beta})$ over the parameter space, which encodes
any knowledge or assumptions about the parameter space that we may
wish to impose before the data are considered, with a likelihood pdf
$\pi_{\text{like}}(\vec{d}_{\text{obs}}|\vec{\beta})$, which explicitly
represents the probability that a given set of parameters $\vec{\beta}$
might give rise to the observed data $\vec{d}_{\text{obs}} \in
\mathbb{R}^m$.  Bayes' Theorem then states the posterior
pdf explicitly as
\begin{align}\label{eq:Bayes}
\pi_{\text{post}}(\vec{\beta} | \data) \propto
\pi_{\text{prior}}(\vec{\beta}) \pi_{\text{like}}(\data | \vec{\beta}).
\end{align}
Note that the infinite-dimensional analog of 
\eqref{eq:Bayes} cannot be
stated 
using pdfs but requires Radon-Nikodym derivatives
\cite{Stuart10}.

For many problems, it is reasonable to choose the prior distribution
to be Gaussian. If the parameters represent a spatial discretization
of a field, the prior covariance operator usually imposes smoothness
on the parameters. This is because rough components of the parameter field
are typically not observable from the data and must be determined by
the prior to result in a well-posed Bayesian inverse problem. Here, we use
elliptic PDE operators to construct the prior covariance, which allows
us to build on fast, optimal complexity solvers.  More
precisely, the prior covariance operator is the inverse of the square
of the Laplacian-like operator~\eqref{eq:reg}, namely $\Cprior :=
\mathcal{A}^{-2} = (-\gamma \gbf
  \Delta_{\Gamma} + \delta I)^{-2}$, where $\gamma$, $\delta > 0$ control the
correlation length and the variance of the prior operator.
This choice of prior ensures that $\Cprior$
is a trace-class operator, guaranteeing bounded pointwise variance and
a well-posed infinite-dimensional Bayesian inverse problem
\cite{Stuart10, Bui-ThanhGhattasMartinEtAl13}.

The difference between the observables predicted by the
model and the actual observations $\data$ (velocity
observations on the ice sheet's surface) is due to
both measurement and model errors, and is represented by the
independent and identically distributed Gaussian random variable
``noise'' vector
$\vec{e} = \vec f(\vec{\beta}) - \vec{d}_{\text{obs}}  \in \mathbb{R}^m,$
where $\vec f(\cdot)$ is the (generally nonlinear) operator mapping
model parameters to output observables. Evaluation of the {\em
  parameter-to-observable map} $\vec f(\vec{\beta})$ requires solution
of the nonlinear Stokes ice flow equations \eqref{eq:forward} followed
by extraction of the surface ice flow velocity field, i.e., it is a
discretization of $\mathcal{B}\uu(\beta)$ in \eqref{objfunction}.
Restating Bayes' theorem with Gaussian noise and prior, we obtain 
the statistical solution of the inverse problem,
$\pi_{\text{post}}(\vec{\beta})$, as
\begin{equation}
\pi_{\text{post}}(\vec{\beta})
\propto\exp\Big(
-\frac{1}{2}\|\vec f(\vec{\beta}) - \vec{d}_{\text{obs}} 
\|^2_{\matrix{\Gamma}_{\text{noise}}^{-1}} 
-\frac{1}{2}\|\vec{\beta}-\mpr\|^2_{\matrix{\Gamma}_{\text{prior}}^{-1}}
\Big),\label{posterior}
\end{equation}
where $\mpr$ is the mean of
the prior distribution, 
$\matrix{\Gamma}_{\text{prior}}\in
\mathbb{R}^{n\times n}$ is the covariance matrix for the prior that
arises upon discretization of $\Cprior$, and
$\matrix{\Gamma}_{\text{noise}}\in \mathbb{R}^{m\times m}$ is the
covariance matrix for the noise, which takes over the role of the
scaling matrix $\matrix\Gamma$ from the deterministic formulation
presented in section~\ref{ipmethod}.

As is clear from the expression \eqref{posterior}, despite the choice
of Gaussian prior and noise probability distributions, the posterior
probability distribution need not be Gaussian, due to the nonlinearity
of $\vec f(\vec{\beta})$. The non-Gaussianity of the posterior poses challenges
for computing statistics of interest for typical large-scale inverse
problems, since
$\pi_{\text{post}}$ is often a surface in high dimensions, and evaluating each
point on this surface requires the solution of the forward model.
Numerical quadrature to compute the mean and covariance
matrix, for example, is completely out of the question. The method of
choice 
is Markov chain Monte Carlo (MCMC), which judiciously samples the
posterior
so that sample statistics can be computed. But the use of MCMC for
large-scale inverse problems is still prohibitive for expensive
forward problems (such as those governed by the nonlinear PDEs of ice
sheet flow) and high dimensional parameter spaces (such as the
$O(10^5-10^6)$ parameters characterizing the basal sliding parameter
field), since even for modest numbers of parameters, the number of
samples required can be in the millions (see for example the
discussion in \cite{PetraMartinStadlerEtAl14,
  FoxHaarioChristen12}). The need to execute millions of forward ice
sheet model simulations is simply not feasible, even with today's
multi-petaflops systems.

We are thus led to make a quadratic approximation of the negative log
of the posterior \eqref{posterior}, which results in a Gaussian
approximation of the posterior $\mathcal{N}(\vec{\map},\postcov)$.
The mean of this posterior distribution, $\vec{\map}$, is the
parameter vector maximizing the posterior \eqref{posterior}, and
is known as the {\em maximum a posteriori} (MAP) point.  It can be found
by
minimizing the negative log of \eqref{posterior}, which amounts to
solving the optimization problem~\eqref{objfunction} (i.e., the
deterministic inverse problem) with appropriately weighted norms,
i.e., 
\begin{equation}\label{eq:objfunction-bayesian}
\vec{\map} := \underset{\vec{\beta}}{\arg \min}
\Big(
-\frac{1}{2}\| \vec f(\vec{\beta}) - \vec{d}_{\text{obs}}
\|^2_{\matrix{\Gamma}_{\text{noise}}^{-1}} 
-\frac{1}{2}\|\vec{\beta}-\mpr\|^2_{\matrix{\Gamma}_{\text{prior}}^{-1}} \Big).
\end{equation}
The
posterior covariance matrix $\postcov$ is then given by the inverse of
the Hessian matrix of $\J$ at $\vec{\map}$, namely
\begin{equation}
\postcov = \left(\Hmisfit(\vec{\map}) + \prcov^{-1} \right)^{-1},
\label{eq:posterior}
\end{equation}
where the Hessian of the misfit is given by
\begin{equation}
\Hmisfit := \FF^*\ncov^{-1}\FF-\FF^*\W_{\!u\beta}-\W_{\!\beta u}\FF + \W_{\!\beta\beta}.
\label{eq:Hmisfit}
\end{equation}
Here $\FF$ is the Jacobian matrix of the parameter-to-observable map
evaluated at $\dmap$, $\FF^*$ is its adjoint, and $\W_{\!u\beta}$,
$\W_{\!\beta u}$, and $\W_{\!\beta\beta}$ involve second derivatives
of the negative log posterior with respect to the parameters and the
states (explicit expressions for the ice sheet inverse problem will be
given below). The Gaussian approximation will be accurate when the
parameter-to-observable map $\vec f(\vec{\beta})$ behaves nearly
linearly over the support of the posterior. This will be the case not
only when $\vec f(\vec{\beta})$ is weakly nonlinear, but also for
directions in parameter space that are poorly informed by the data (in
which case $\vec f(\vec{\beta})$ is approximately constant and thus
the prior dominates), as well as directions in parameter space that
are strongly informed by the data (in which case the posterior
variance is small and thus the linearization of $\vec f(\vec{\beta})$
is accurate over the support of the posterior). For the present work,
in which we tackle a massive scale Bayesian inverse problem ($\sim
10^6$ parameters) governed by large-scale ice sheet flow in
Antarctica, all known MCMC methods will be prohibitive, so we cannot
compare our Hessian-at-the-MAP-based Gaussian approximation with full
sampling of the posterior. However, in
\cite{PetraMartinStadlerEtAl14}, we compare this Gaussian
approximation with several variants of MCMC sampling of the posterior
for a model 2D ice flow problem with a moderate number of parameters
($\sim$100), and conclude that the Hessian-based Gaussian
approximation can be appropriate, both as a proposal for MCMC, as well
as a stand-alone approximation of the posterior.

\subsection{Low-rank based posterior covariance approximation}
\label{sec:low-rank}

\begin{figure}
\begin{center}
\includegraphics[width=.4\columnwidth]{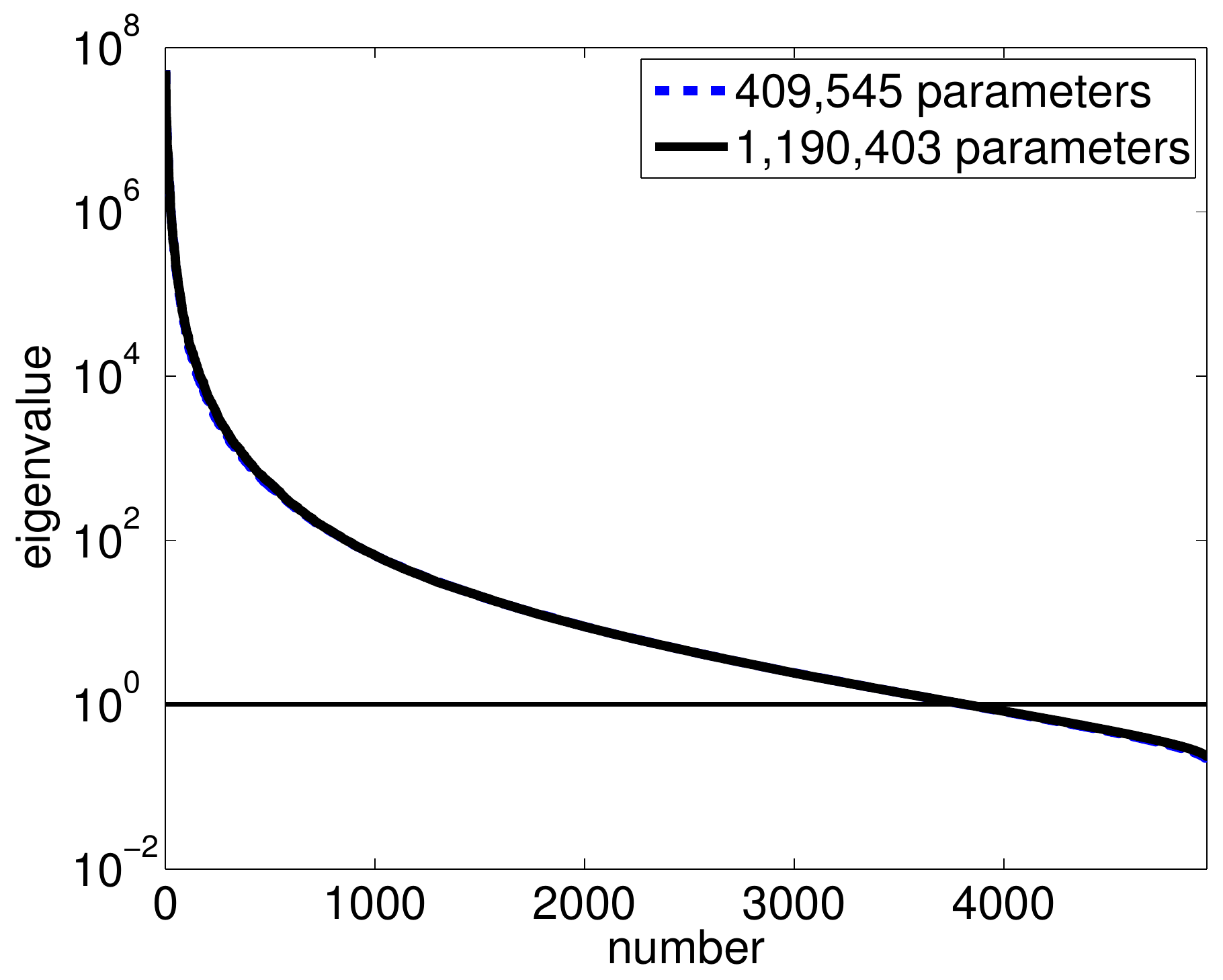}
\end{center}
\caption{Log-linear plot of spectrum of prior-preconditioned data
  misfit Hessian for two successively finer parameter/state meshes of
  the inverse ice sheet
  problem~\eqref{eq:objfunction-bayesian}.
  The curves lie on top of each other, indicating mesh
  independence of the spectrum. This implies that the dominant
  eigenvalues/eigenvector of the parameter field---and thus the
  information content of the data, filtered through the prior---are
  independent of the parameter dimension. Moreover, since the surface
  observations are refined with the state mesh, the invariance of the
  spectrum to mesh refinement also implies that the dominant
  eigenvalues are independent of the data dimension (once the
  information content of the data is resolved). The low rank
  approximation captures this dominant, data-informed portion of the
  spectrum. The eigenvalues are truncated at 0.2.
}
\label{fig:spectrum}
\end{figure}

As stated above, the Gaussian approximation of the posterior
\eqref{posterior}, with covariance matrix $\postcov$
\eqref{eq:posterior} that involves the Hessian of the data misfit
evaluated at the MAP point \eqref{eq:Hmisfit}, is still intractable.
The primary difficulty is that the large parameter dimension $n$
prevents any representation of the posterior covariance $\postcov$ as
a dense operator. In particular, the Jacobian of the
parameter-to-observable map, $\FF$, is formally a dense matrix, and
requires $n$ forward PDE solves to construct explicitly. This is
intractable when $n$ is large and the forward PDEs are expensive, as
in our case of Antarctic ice sheet flow. However, a key feature of the
operator $\FF$ is that its action on a (parameter field-like) vector
can be formed by solving a (linearized) forward PDE problem;
similarly, the action of its adjoint $\FF^*$ on a (observation-like)
vector can be formed by solving a (linearized) adjoint PDE.
Moreover, for many ill-posed inverse problems, the Hessian matrix of
the data misfit term in \eqref{objfunction}, given
by~\eqref{eq:Hmisfit}, is a discretization of a compact operator,
i.e., its eigenvalues collapse to zero.  This can be understood
intuitively, since only the modes of the parameter field that strongly
influence the ice velocity are present in the dominant spectrum of
\eqref{eq:Hmisfit}.  In many ill-posed inverse problems, numerous
modes of the parameter field (for example, highly oscillatory ones)
will have negligible effect on the observables. The range space thus
is effectively finite-dimensional even before discretization (and
therefore independent of any mesh), and the eigenvalues decay, often
rapidly, to zero.  Figure~\ref{fig:spectrum} illustrates the rapid
spectral decay of the data misfit Hessian stemming from our Antarctic
ice sheet inverse problem, demonstrating that
only about 5,000 (out of
1,190,403) modes of the basal sliding parameter field can
be inferred from the data (on the fine mesh). Next we exploit this
low-rank structure and the ability to form matrix-free Hessian actions
to construct scalable algorithms to approximate the posterior
covariance operator.

Rearranging the expression for $\postcov$ in \eqref{eq:posterior} to
factor out $\prcov$ gives
\begin{align}
  \label{eqn:gamma_post}
  \postcov = \left(\prcov \Hmisfit +
  \matrix{I}\right)^{-1}\prcov.
\end{align}
This factorization exposes the {\em prior-preconditioned Hessian of
  the data misfit},
\begin{equation}\label{eq:ppmisfitH}
  \HT := \prcov \Hmisfit.
\end{equation}
Note that $\HT$ is symmetric with respect to the inverse
prior-weighted inner product. If a square root of $\prcov$ is
available, one can alternatively use a symmetric preconditioning of
the data misfit Hessian by the prior covariance operator
\cite{FlathWilcoxAkcelikEtAl11, Bui-ThanhGhattasMartinEtAl13}. This
results in a prior-preconditioned data misfit Hessian that is symmetric
with respect to the Euclidean inner product.

To construct a low-rank approximation of $\HT$, we consider the
generalized eigenvalue problem 
\begin{equation}
 \Hmisfit \matrix{W} = \prcov^{-1} \matrix{W} \matrix{\Lambda},
  \label{eq:gevd}
\end{equation}
where $\matrix{\Lambda} = \diag(\lambda_i) \in \mathbb{R}^{n\times n}$
contains the generalized eigenvalues and the columns of $\matrix W\in
\mathbb R^{n\times n}$ the generalized eigenvectors. Because $\prcov$
is symmetric positive definite and $\Hmisfit$ is symmetric positive
semi-definite, $\matrix{\Lambda}$ is non-negative and $\matrix{W}$ is
orthogonal with respect to $\prcov^{-1}$, so we can scale $\matrix W$ so
that $\matrix W^T \prcov^{-1} \matrix W = \matrix I$.  Defining
$\matrix V := \prcov^{-1}\matrix W$, we can rearrange \eqref{eq:gevd}
to get
\begin{equation}
  \prcov \Hmisfit = \matrix{W} \matrix{\Lambda} \matrix{V}^T.
  \label{eq:HT}
\end{equation}
When the generalized eigenvalues $\{\lambda_i\}$ decay rapidly, we can
extract a low-rank approximation of $\HT$ by retaining only the $r$
largest eigenvalues and corresponding eigenvectors,
\begin{equation}
\label{eq:HT_approx}
  \HT \approx
  \matrix{W}_r \matrix{\Lambda}_r \matrix{V}_r^T.
\end{equation}
Here, $\matrix{W}_r \in \mathbb{R}^{n\times r}$ contains only the $r$
eigenvectors of $\HT$ that correspond to the $r$ largest eigenvalues,
which are assembled into the diagonal matrix $\matrix{\Lambda}_r = \diag
(\lambda_i) \in \mathbb{R}^{r \times r}$, and
$\matrix{V}_r=\prcov^{-1} \matrix{W}_r$.

Here, we use randomized SVD algorithms \cite{HalkoMartinssonTropp11,
  LibertyWoolfeMartinssonEtAl07} to construct the approximate spectral
decomposition. Their performance is comparable to Krylov methods such
as Lanczos, which has been used for low-rank approximations of
Hessians in very high dimensions in \cite{FlathWilcoxAkcelikEtAl11,
  Bui-ThanhGhattasMartinEtAl13, KalmikovHeimbach14}.
Both algorithms require only the action of the Hessian on
vectors; we show how to do this for the ice sheet Hessian operator in
Section~\ref{hess}.  However, randomized methods have a significant
edge over deterministic methods for large-scale problems, since the
required Hessian matrix-vector products are independent of each other,
providing asynchronicity and fault tolerance.

\begin{figure}[ht]
\centering
\begin{minipage}{0.99\textwidth}
  \includegraphics[width=.32\columnwidth]{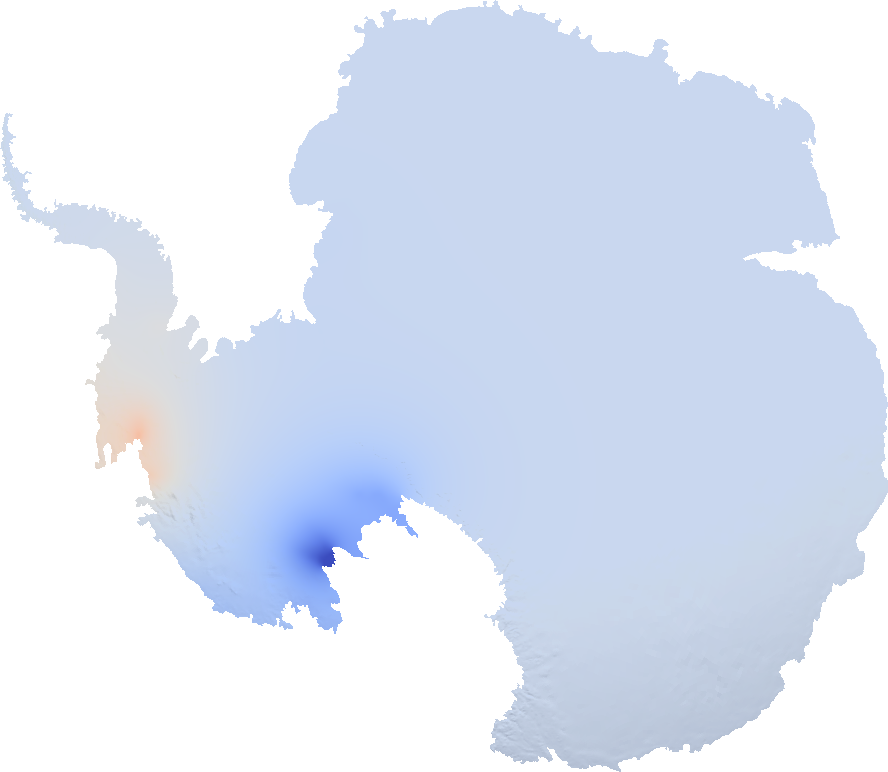} \hfill
  \includegraphics[width=.32\columnwidth]{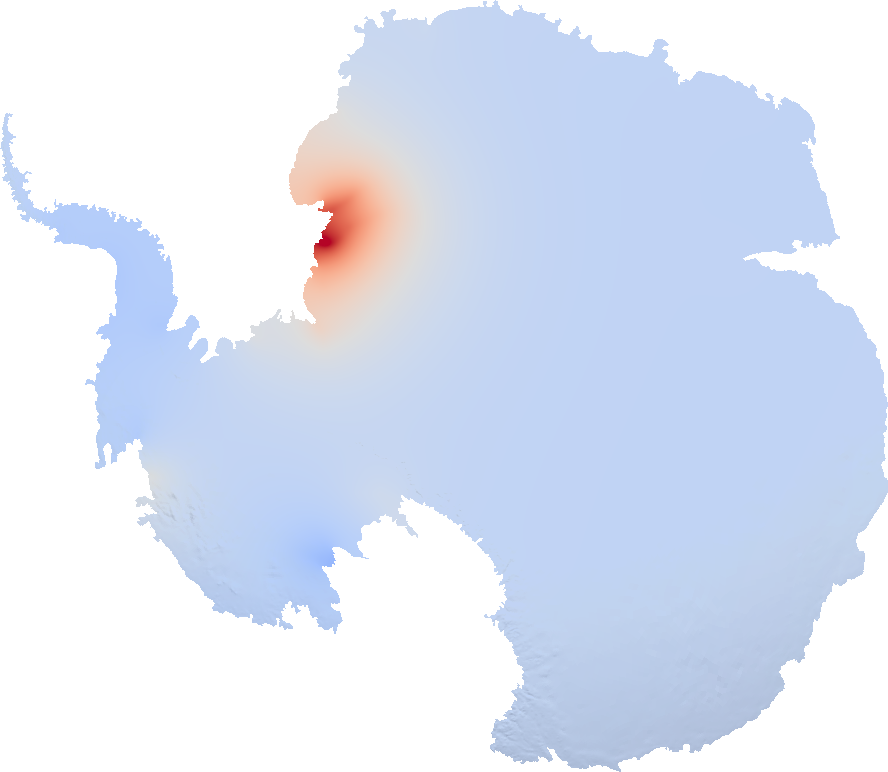} \hfill
  \includegraphics[width=.32\columnwidth]{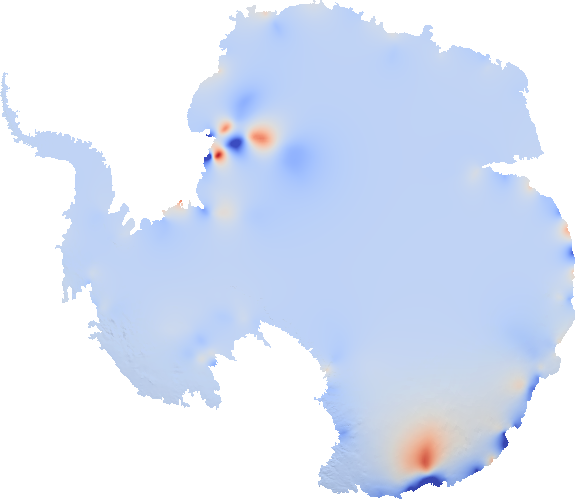} \\
  \includegraphics[width=.32\columnwidth]{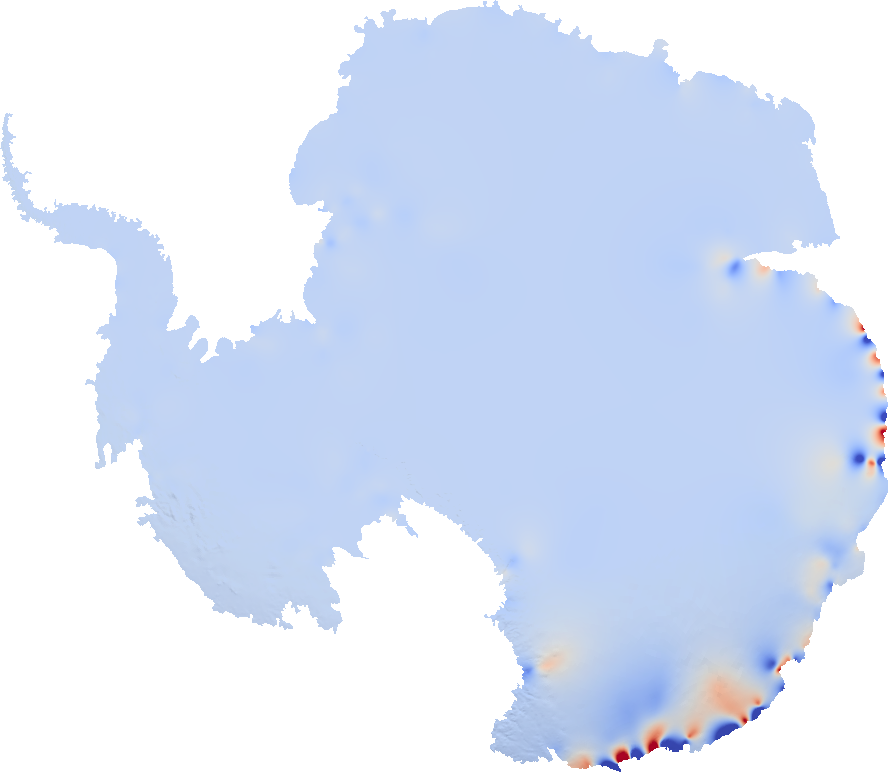} \hfill
  \includegraphics[width=.32\columnwidth]{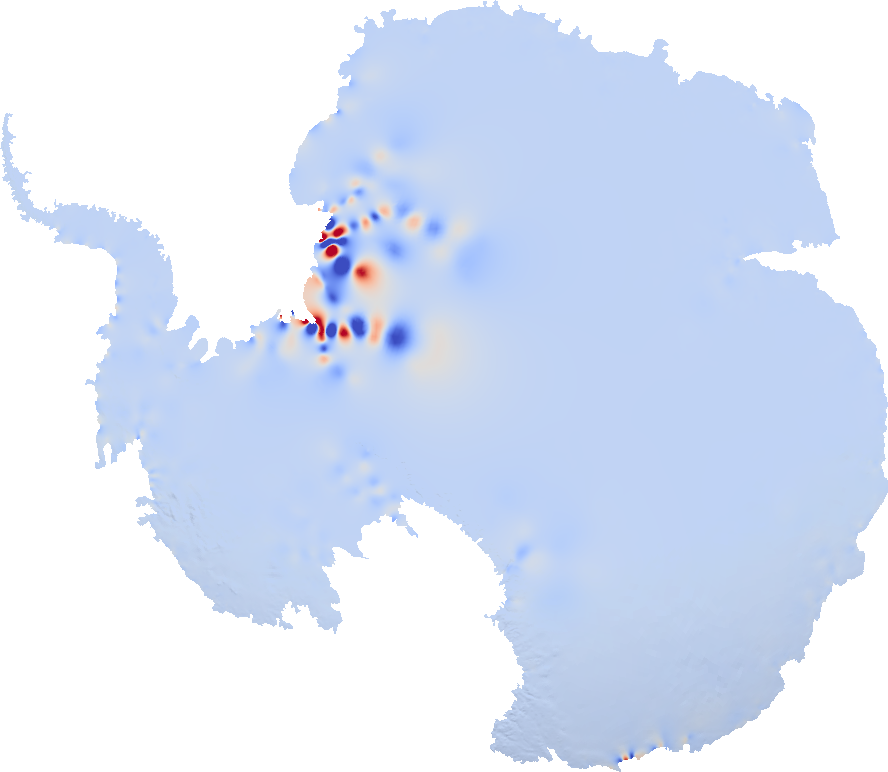} \hfill
  \includegraphics[width=.32\columnwidth]{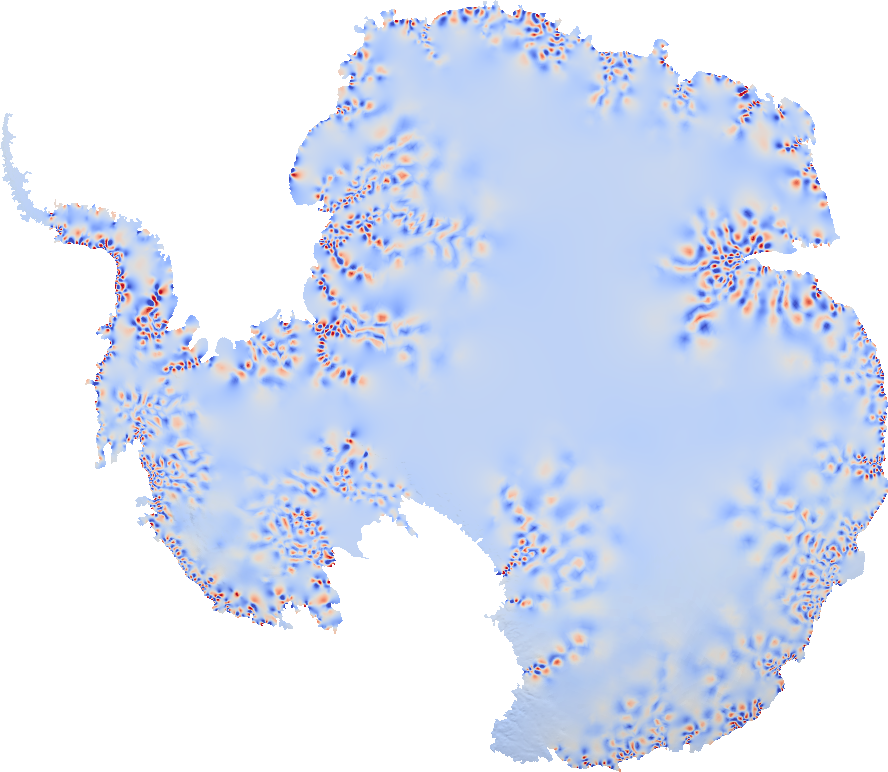} \hfill
\end{minipage}
\caption{Eigenvectors of the prior-preconditioned data misfit Hessian
  corresponding (from left to right and top to bottom) to the 1st, 2nd, 100th, 200th, 500th, and 4000th
  eigenvalues. Note that different eigenvectors localize in different
  regions of Antarctica, and that eigenvectors corresponding to
  smaller eigenvalues are increasingly more oscillatory (and thus
  inform smaller length scales of the basal sliding parameter field) but
  are also increasingly less informed by the data.}
\label{fig:evecs}
\end{figure}

Once the low-rank approximation \eqref{eq:HT_approx} has been
constructed, we proceed to obtain the posterior covariance matrix.
The Sherman-Morrison-Woodbury formula is employed to perform the
inverse in \eqref{eqn:gamma_post},
\begin{align}
  \label{eqn:eigen_error}
  \notag \left(\HT+ \matrix{I}\right)^{-1}
  = \matrix{I}-\matrix{W}_r \matrix{D}_r \matrix{V}_r^T +
  \mathcal{O}\left(\sum_{i=r+1}^{n} \frac{\lambda_i}{\lambda_i +
    1}\right), \nonumber
\end{align}
where $\matrix{D}_r :=\diag(\lambda_i/(\lambda_i+1)) \in
\mathbb{R}^{r\times r}$. The last term in this expression captures the
error due to truncation in terms of the discarded eigenvalues; this
provides a criterion for truncating the spectrum, namely that $r$ is
chosen such that $\lambda_r$ is small relative to 1. With this
low-rank approximation, the final expression for the approximate
posterior covariance follows from \eqref{eqn:gamma_post},
\begin{equation}
\label{eqn:reduced_post}
\postcov \approx (\matrix{I} - \matrix{W}_r \matrix{D}_r
\matrix{V}_r^T) \prcov = 
\prcov
- \matrix{W}_r \matrix{D}_r \matrix{W}_r^T.
\end{equation}
Note that \eqref{eqn:reduced_post} expresses the posterior uncertainty
as the prior uncertainty less any information gained from the data,
filtered through the prior (as a consequence of choosing the
prior-preconditioned data misfit Hessian as the operator whose
spectrum is truncated). The retained eigenvectors of the
prior-preconditioned data misfit Hessian are those modes in parameter
space that are simultaneously well-informed by the data and assigned
high probability by the prior.  Figure~\ref{fig:evecs} displays
several of these eigenvectors.
The low rank update of the prior covariance in
\eqref{eqn:reduced_post}, which is based on the low rank approximation
\eqref{eq:HT_approx} of the prior-preconditioned Hessian of the data
misfit, has recently been shown to be the optimal rank-$r$ update with respect to several
important loss functions \cite{SpantiniSolonenCuiEtAl14}.

\subsection{The action of the Hessian operator}
\label{hess}

In this section we show that the Hessian-vector products needed in the
randomized low-rank algorithms described above can be computed
efficiently and matrix-free. The action of the Hessian operator
evaluated at a basal sliding parameter field $\beta$ in a direction $\betah$ is
given by
\begin{align}
\label{eq:hessian_action}
  \mc{H} (\beta) \betah  &:= \exp(\beta) (
  \betah \gbf T \uu \cdot \gbf T \vv +
  \gbf T \uh \cdot \gbf T \vv +
  \gbf T \uu \cdot \gbf T \vh) +
  \mc{A}^{2\kappa}(\betah),
\end{align}
where $(\uh,\ph)$ satisfy a certain linearized forward Stokes equation, which
we call the {\it incremental forward Stokes
  equations}
\begin{subequations}\label{eq:incforward}
\begin{alignat}{2}
  \gbf{\nabla} \cdot  \uh & = \,0 & &  \quad \text{ in } \Omega\\
  - \gbf{\nabla} \cdot \gbf{\sigma}_{\!\uh} &= \gbf{0} &
  &  \quad \text{ in } \Omega\\
  \gbf{\sigma}_{\!\uh} \gbf{n} & = \gbf{0} &
  &  \quad \text{ on }\Gamma_{\!\mbox{\tiny} t}\\
  \uh \cdot \gbf{n} &=  0 & & \quad  \text{ on } \Gamma_{\!\mbox{\tiny} b} \\
  \gbf T\gbf{\sigma}_{\!\uh} \gbf{n} +\exp(\beta) \gbf T
  \uh  & = - \betah  \exp(\beta) \gbf T \uu & & \quad
  \text{ on } \Gamma_{\!\mbox{\tiny} b},
\end{alignat}
\end{subequations}
with the incremental forward stress
\begin{align*}
  \gbf{\sigma}_{\!\uh} := 2
  \eta(\uu)\, \bigl(\mathsf{I} + \frac{1-n}{n} \frac{\edot_{\uu}\otimes
    \edot_{\uu}}{\edot_{\uu} \cdot \, \edot_{\uu}}\bigr)\edot_{\uh
  }-\gbf{I} \ph,
\end{align*}
and $(\vh, \qh)$ satisfy a linearized version of the adjoint equation,
the {\it incremental adjoint Stokes equations}
\begin{subequations}\label{eq:incadjoint}
\begin{alignat}{2}
  \gbf{\nabla} \cdot \vh & = \,0  &\:
  &\text{ in } \Omega\\
  - \gbf{\nabla} \cdot \gbf{\sigma}_{\!\vh} &=
  - \gbf{\nabla} \cdot \gbf{\tau}_{\!\uh} & &\text{ in }
  \Omega\\
  \gbf{\sigma}_{\!\vh} \gbf{n} & =
  -\mathcal{B}^*\bs\Gamma\mathcal{B}\uh
  - \gbf{\tau}_{\!\uh}\gbf{n} &\: &\text{ on }
  \Gamma_{\!\mbox{\tiny} t}\\
  \vh\cdot \gbf{n} &= 0 & & \text{ on }
  \Gamma_{\!\mbox{\tiny} b}\\
  \gbf T\gbf{\sigma}_{\!\vh}
  \gbf{n} + \exp(\beta) \gbf T \vh &= - \gbf T
  \gbf{\tau}_{\!\uh}\gbf{n} & & \text{ on }
  \Gamma_{\!\mbox{\tiny} b},
\end{alignat}
where the incremental adjoint stress $\gbf{\sigma}_{\!\vh}$ is given
by
\begin{align*}
\gbf{\sigma}_{\!\vh} := 2 \eta(\uu) \,
\bigl(\mathsf{I} + \frac{1-n}{n} \frac{\edot_{\uu}\otimes
  \edot_{\uu}}{\edot_{\uu} \cdot \,
  \edot_{\uu}}\bigr)\edot_{\vh} -\gbf{I}\qh,
\end{align*}
and $\gbf{\tau}_{\!\uh} = 2 \eta(\uu) \Psi \edot_{\uh}$, where
\begin{align*}
  \Psi & := \frac{1-n}{n} \biggr[ \frac{\edot_{\uu}\otimes
      \edot_{\vv} + \edot_{\vv}\otimes \edot_{\uu}}{\edot_{\uu} \cdot \,\edot_{\uu}} +
    \frac{\edot_{\uu}\cdot \, \edot_{\vv}}{\edot_{\uu} \cdot \,\edot_{\uu}}
    \biggl( \frac{1-3n}{n} \frac{\edot_{\uu}\otimes \edot_{\uu}}{\edot_{\uu} \cdot \,\edot_{\uu}} + \mathsf{I} \biggr)\biggr].
\end{align*}
\end{subequations}
We see that the incremental forward and incremental adjoint Stokes
equations are linearized versions of their forward and adjoint
counterparts, having the same (linearized) operator and differing only
in the source terms. Thus, computation of Hessian actions
on vectors amount to solution of a pair of forward/adjoint
(linearized) Stokes equations, similar to the computation of the
gradient. Since the gradient and Hessian action require just one and
two linearized Stokes solves, respectively, they are significantly
cheaper to evaluate than solving the (nonlinear) forward problem,
which typically requires an order of magnitude more linearized Stokes
solves.

\subsection{Scalability of Bayesian solution of the inverse problem}
\label{scalability}

We now discuss the overall scalability of our algorithms for Bayesian
solution of the inverse problem.
First, we characterize the scalability of construction of the
low-rank-based approximate posterior covariance matrix in
\eqref{eqn:reduced_post}. As stated before, the
parameter-to-observable map $\FF$ cannot be constructed explicitly,
since it requires $n$ linearized forward Stokes solves. However, as
elaborated above, its action on a vector can be computed by a single
linearized Stokes solve, regardless of the number of parameters $n$
and observations $m$. Similarly, the action of $\FF^*$ on a vector can
be computed via a linearized adjoint Stokes solve.  Moreover, the
prior is usually much cheaper to apply than a forward/adjoint solve
(here, it is a scalar elliptic solve on the basal boundary).
Therefore, the cost of applying $\HT$ to a vector---and thus the per
iteration cost of the randomized SVD algorithm---is dominated by a
pair of linearized forward and adjoint Stokes solves.

The randomized SVD algorithm requires a number of matrix-vector
products proportional to the effective rank $r$ of the matrix
\cite{LibertyWoolfeMartinssonEtAl07, HalkoMartinssonTropp11}.  Thus,
the remaining component to establish scalability of the low-rank
approximation of $\HT$ is independence of $r$---and therefore the
number of matrix-vector products, and hence Stokes solves---from the
parameter dimension $n$.  This is the case when
$\matrix{H}_{\text{misfit}}$ in \eqref{eq:Hmisfit} is a
(discretization of a) compact operator, and when preconditioning by
$\prcov$ does not destroy the spectral decay. This situation is
typical for many ill-posed inverse problems, in which the prior is
either neutral or of smoothing type.
Hence, a low-rank approximation of $\HT$ can be made that does not
scale with parameter or data dimension, instead depending only on the
information content of the data, filtered through the prior. This is
indeed the case for our ice sheet inverse problem, as demonstrated in 
Figure~\ref{fig:spectrum}. 

Once the $r$ eigenpairs defining the low rank approximation have been
computed, estimates of uncertainty can be computed by interrogating
$\postcov$ in \eqref{eqn:reduced_post} at a cost of just $r$ inner
products (which are negligible) plus elliptic solves representing the
action of the prior $\prcov$ on a vector (here carried out with an
algebraic multigrid solver and therefore algorithmically
scalable). For example, samples can be drawn from the Gaussian defined
with a covariance $\postcov$, a row/column of $\postcov$ can be
computed, and the action of $\postcov$ in a given direction can be
formed, all at cost that is $O(rn)$ in the number of inner products in
addition to the $O(n)$ cost of the multigrid solve. Moreover, the
posterior variance field, i.e., the diagonal of $\postcov$, can be
found with $O(rn)$ linear algebra plus $O(r)$ multigrid solves.

\subsection{Samples from the prior and Gassianized posterior}
\label{eq:uqres}

In this section we provide some results that illustrate the
quantification of uncertainty in the solution of the Antarctic ice
sheet inverse problem. The number of degrees of freedom was 3,785,889
for each of the discretized state-like variables (state, adjoint,
incremental state, and incremental adjoint velocities and pressures)
and 409,545 for the uncertain basal sliding parameter field. The
problem was solved on 1024 processor cores.
We used the Laplacian prior operator with $\kappa = 1$, $\gamma = 10$,
and $\delta = 10^{-5}$, and took a zero mean of the basal sliding
parameter field, i.e., $\beta_0 = 0$. To characterize the noise, we
used a diagonal noise covariance matrix $\ncov$ with entries
$\sigma^i_{\text{noise}} = 0.1(\|\data^i\|^2 + \varepsilon)^{1/2}$ for
$i = 1,...,m$, with $\varepsilon = 10^{-9}$. That is,
$\sigma^i_{\text{noise}}$ is (at least) 10\% of the velocity magnitude
at the $i$th observational data point.

To compute the MAP point via solution of the optimization problem
\eqref{eq:objfunction-bayesian}, 213 (outer) Newton iterations were
necessary to decrease the nonlinear residual by a factor of $10^3$. In
each of these outer Newton iterations, the nonlinear forward Stokes
problem has to be solved, for which we use an (inner) Newton
method. These inner Newton solves are terminated after the residual is
decreased by a factor of $10^8$, which takes an average of 5 Newton
iterations, each requiring a linearized Stokes solve. In addition to
the nonlinear Stokes solve, each (outer) Newton iteration requires
computation of the gradient, which costs one linearized Stokes solve,
and inexact solution of the (outer) Newton
system~\eqref{eq:newton-opt} using CG, which requires a number of
Hessian-vector products.  The average number of CG iterations per
(outer) Newton step was 239. Each Hessian-vector product requires a
pair of linearized Stokes solves, one incremental forward and one
incremental adjoint. The tolerance for these incremental Stokes solves
was set to $10^{-6}$. 
Altogether, finding the MAP point required a total of 107,578
linear(ized) Stokes solves.

We approximated the covariance matrix at the MAP point via a low-rank
representation employing 5,000 products of the Hessian matrix with
random vectors; hence, the cost of this approximation is 10,000
incremental forward/adjoint (linearized) Stokes solves. Thus, the cost
(measured in Stokes solves) of quantifying the uncertainty in the
MAP solution is about an order of magnitude less than that of finding
the MAP point. 

\begin{figure}[t!]
  \centering
\begin{tikzpicture}
  \node (l) at (0,0.5)
  {\includegraphics[width=.6\columnwidth]{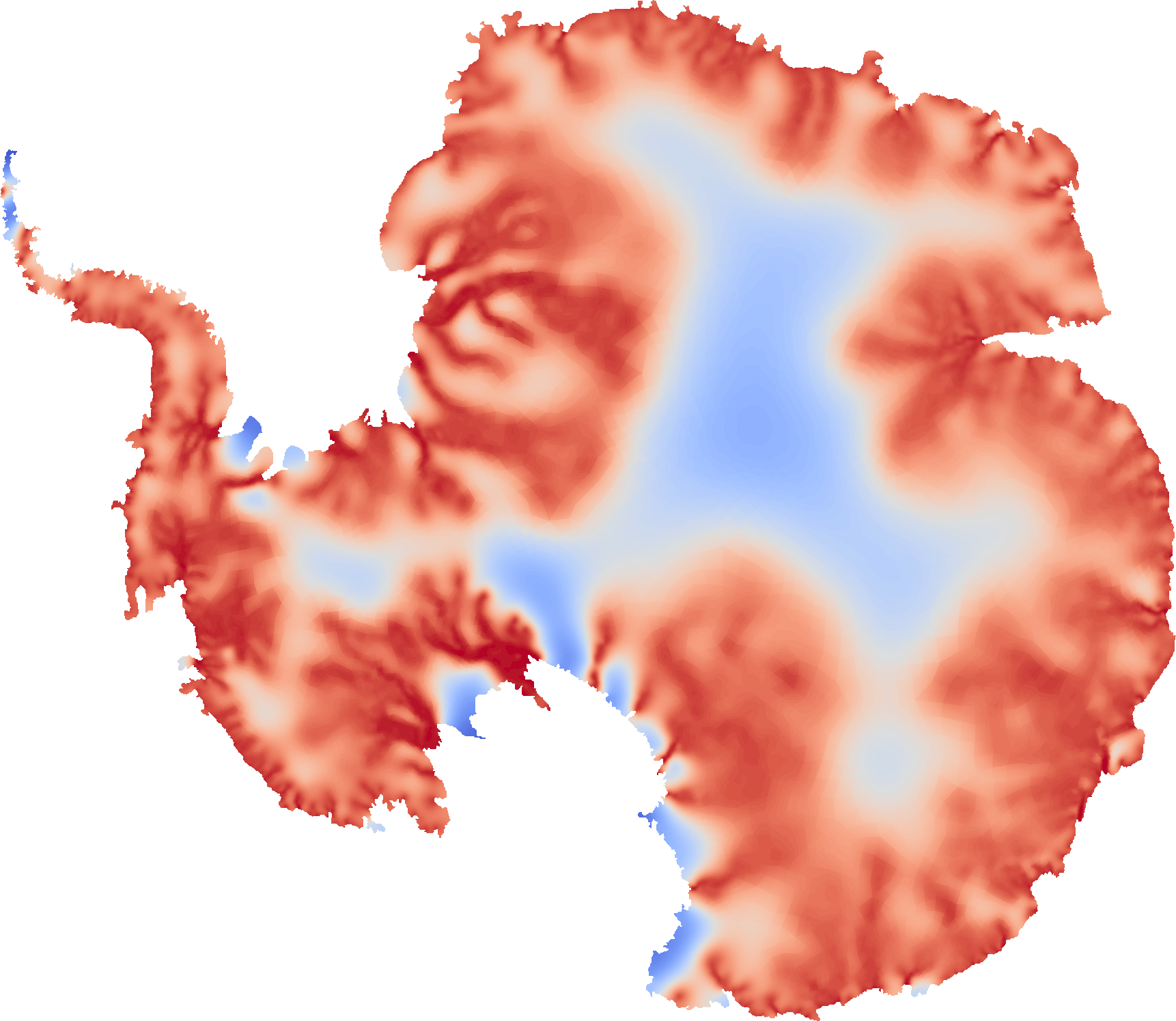}};
  \node (r) at (6,0.5)
  {\includegraphics[width=0.13\columnwidth]{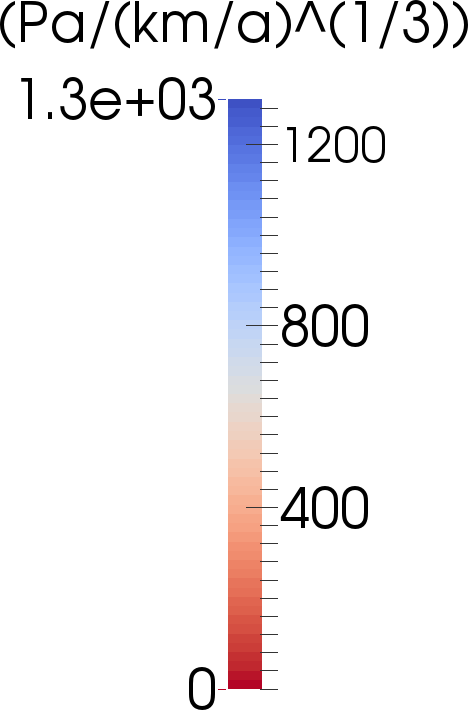}};
  \draw[fill=white,white] (4.8,1.8) rectangle (7.2,2.5);
  \node at (6,2.5) {$\exp(\beta)^{1/3}$};
  \node at (6,2.05) {[(Pa/(km/a))$^{1/3}$]};
\end{tikzpicture}
\caption{The posterior mean. Low (red) and high (blue) values for the basal
  sliding parameter correlate with fast and slow ice flow regions,
  respectively.}
\label{fig:mean}
\end{figure}
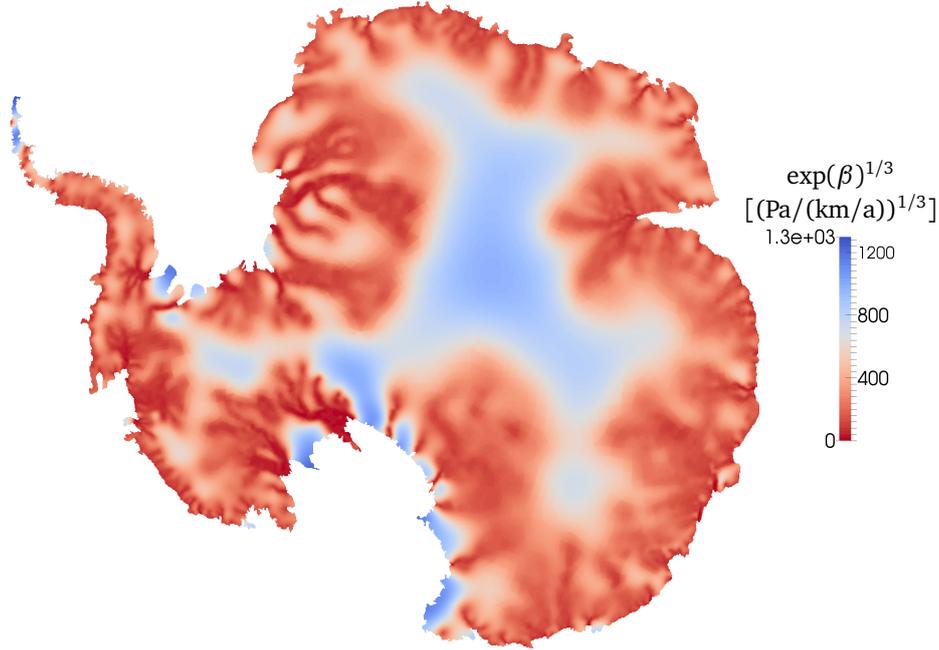
\begin{figure}[ht]
  \centering
  \begin{tikzpicture}
    \node (l) at (0,0.5){
      \includegraphics[width=0.4\columnwidth]{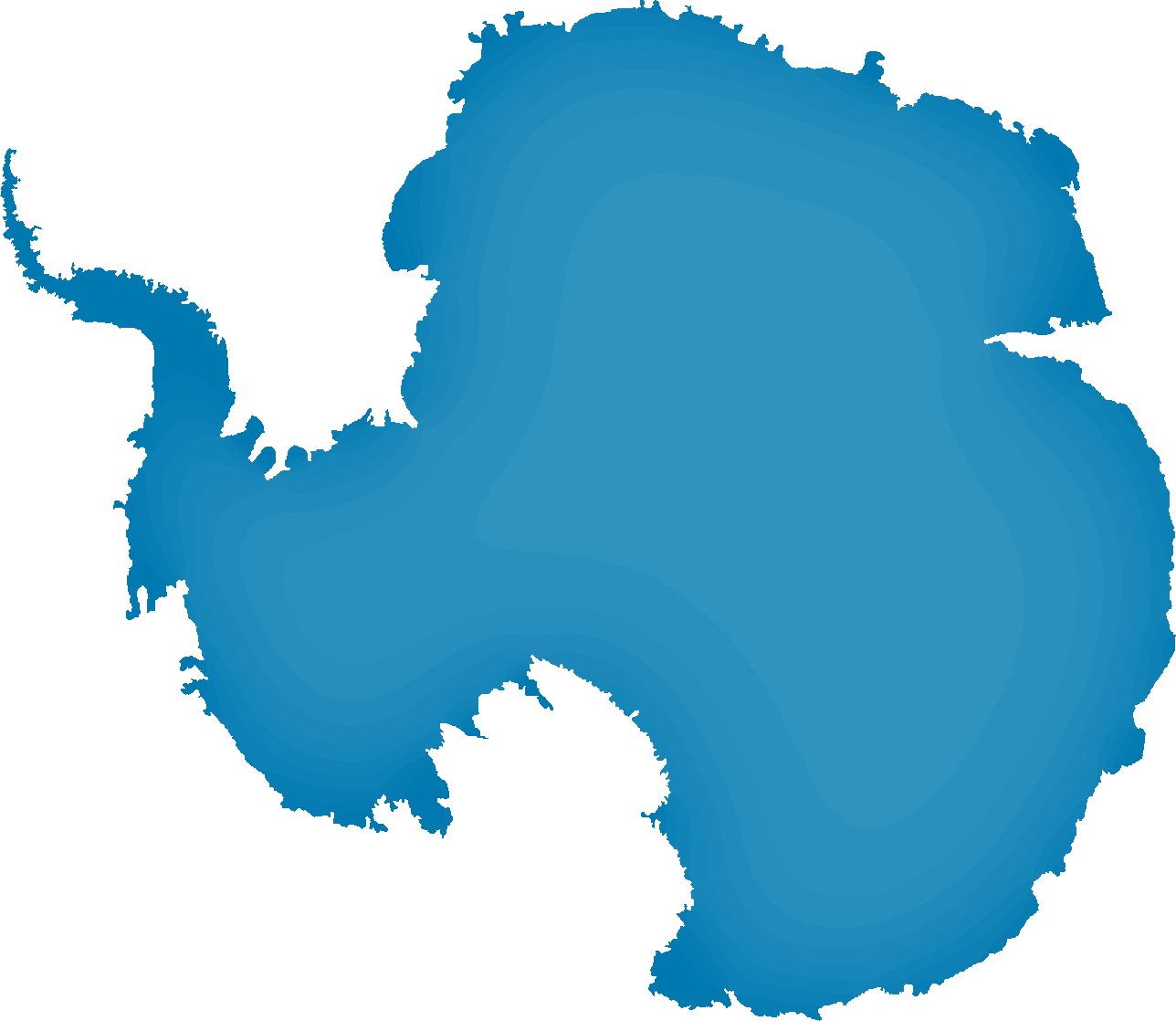}
      \vspace{0.2in}
      \includegraphics[width=0.4\columnwidth]{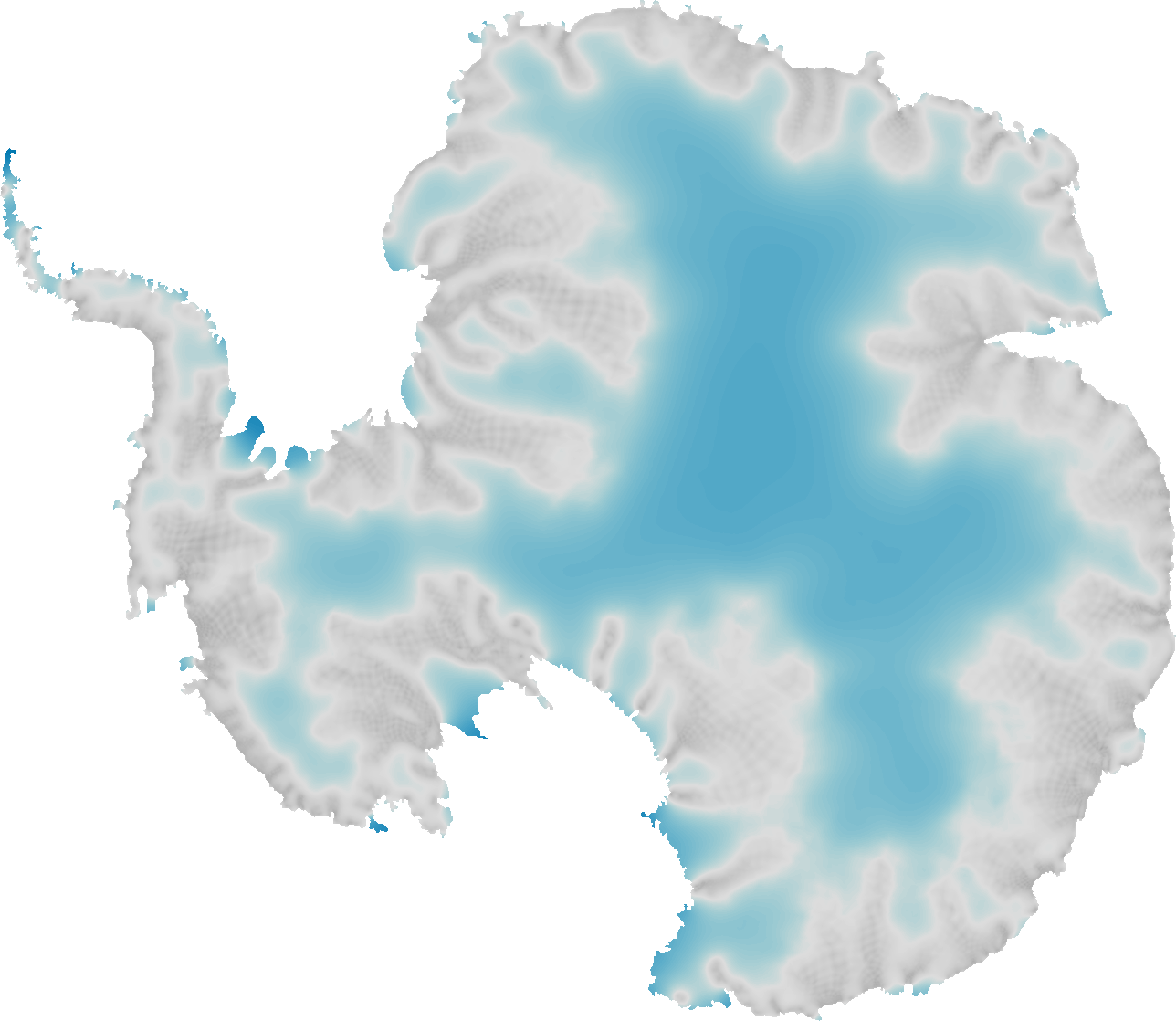}};
    \node (r) at (7.5,0.5){
      \includegraphics[width=0.07\columnwidth]{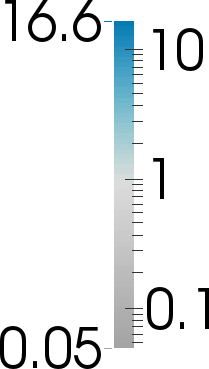}};
    \node at (7.6,1.8) {std. dev. of $\beta$};
  \end{tikzpicture}
\caption{This figure shows the standard deviations of the pointwise
  marginals of the prior distribution (left) and of (the Gaussian
  approximation of) the posterior distribution (right).}
\label{fig:variance}
\end{figure}
\begin{figure}[ht]
  \centering

\begin{tikzpicture}
  \node (l) at (0,0.5){
    \includegraphics[width=.28\columnwidth]{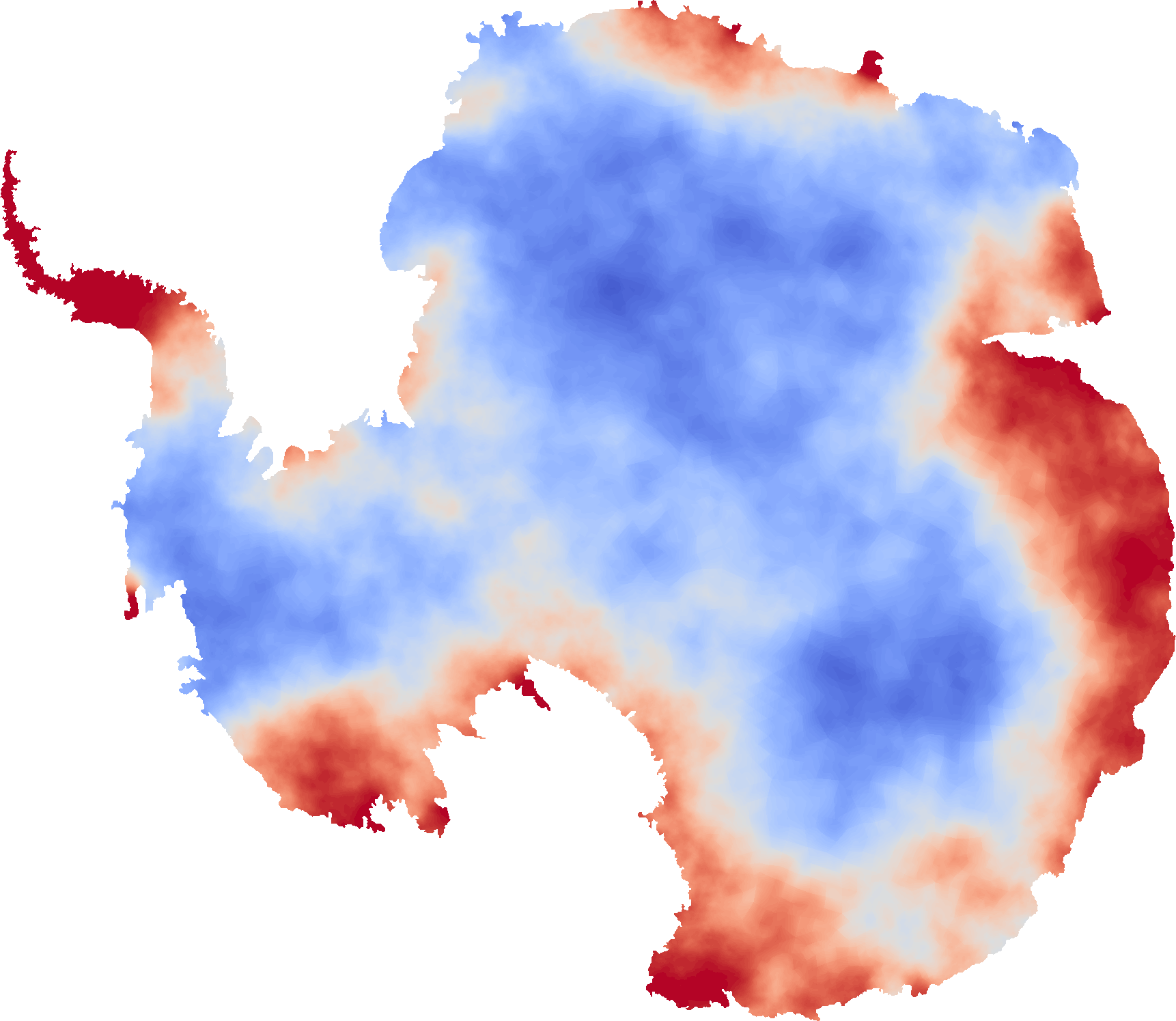}
    \vspace{0.05in} %
    \includegraphics[width=.28\columnwidth]{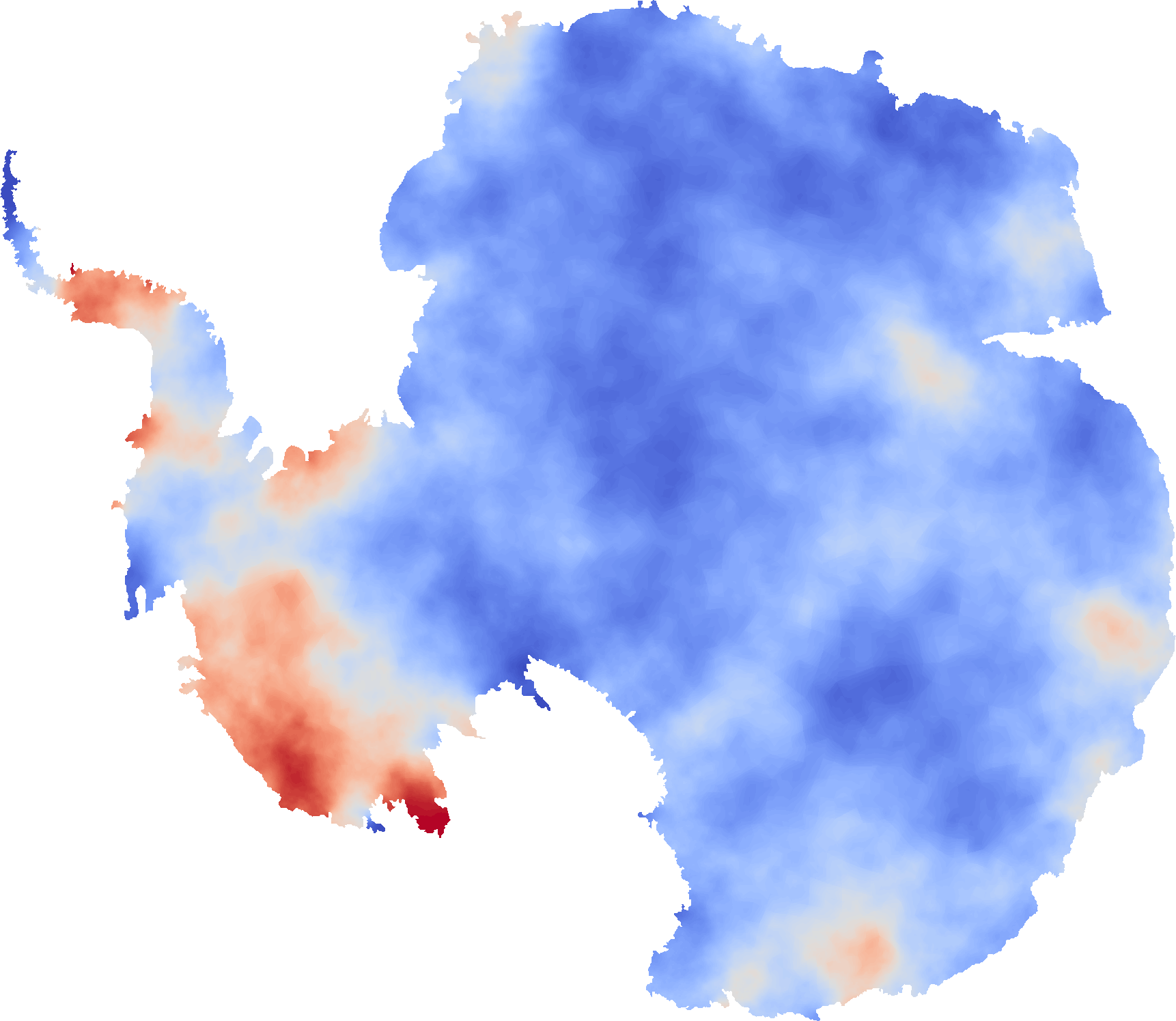}
    \vspace{0.05in} %
    \includegraphics[width=.28\columnwidth]{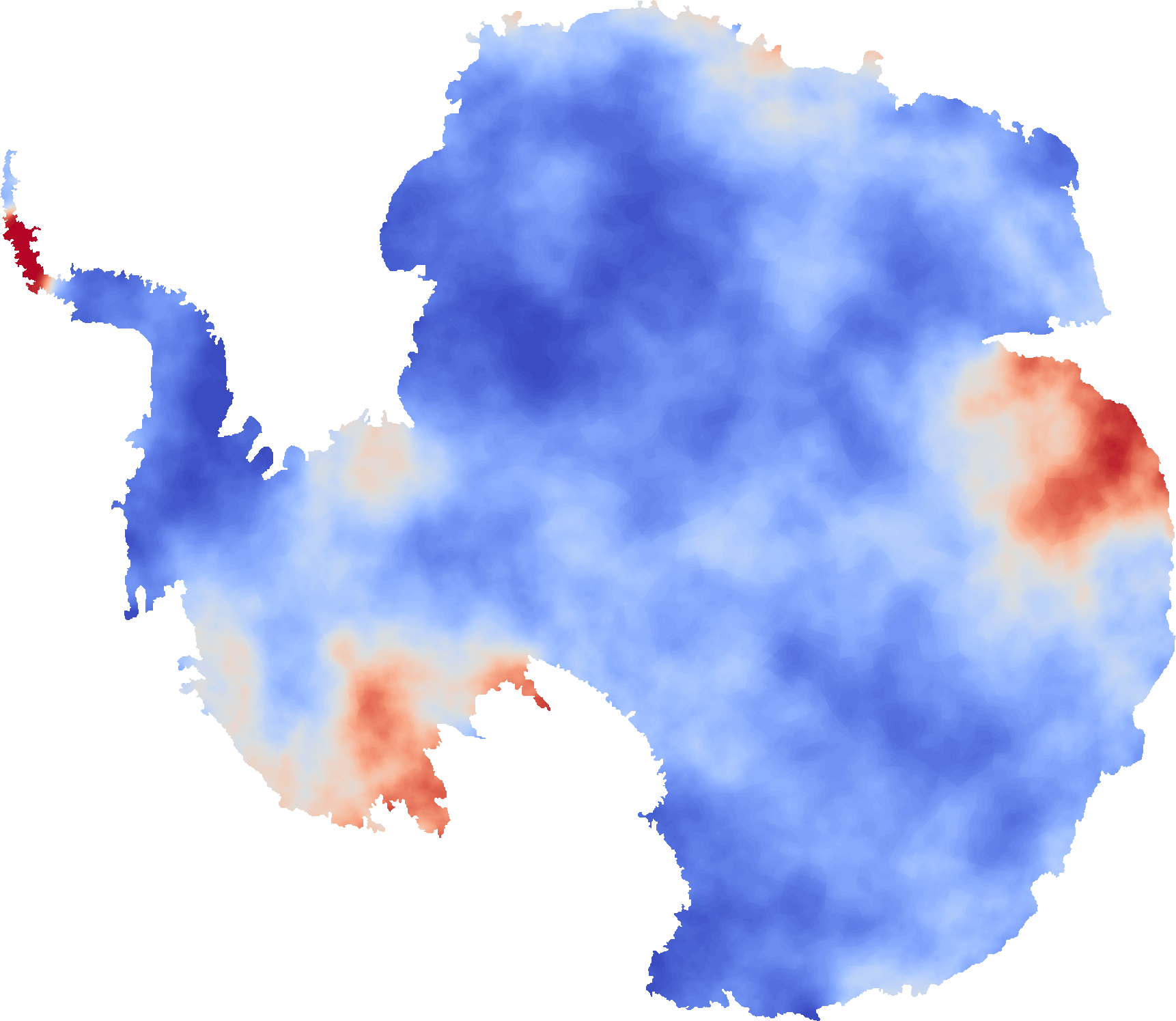}
  };
  \node (ll) at (0,-4){
    \includegraphics[width=.28\columnwidth]{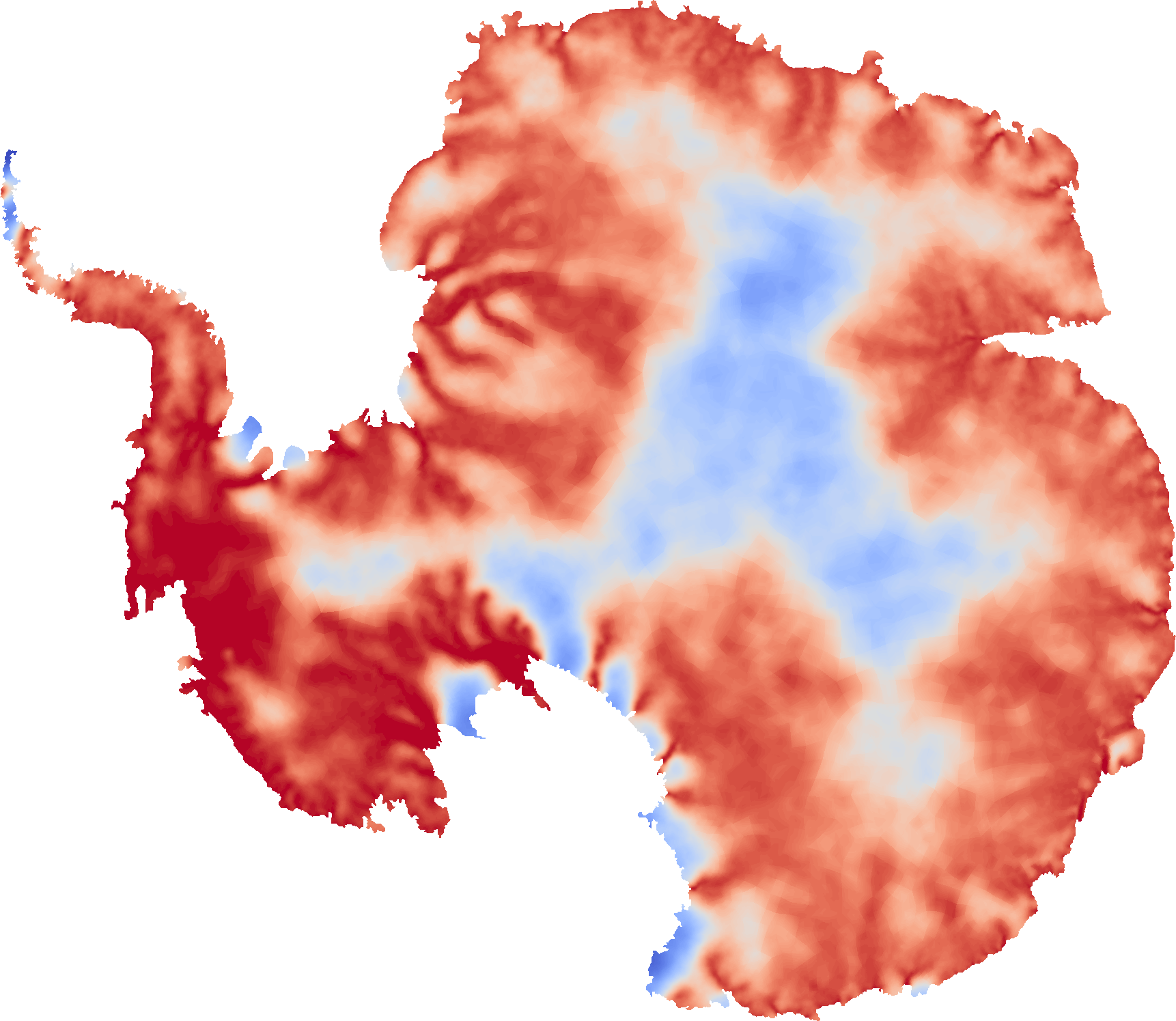}
    \vspace{0.05in} %
    \includegraphics[width=.28\columnwidth]{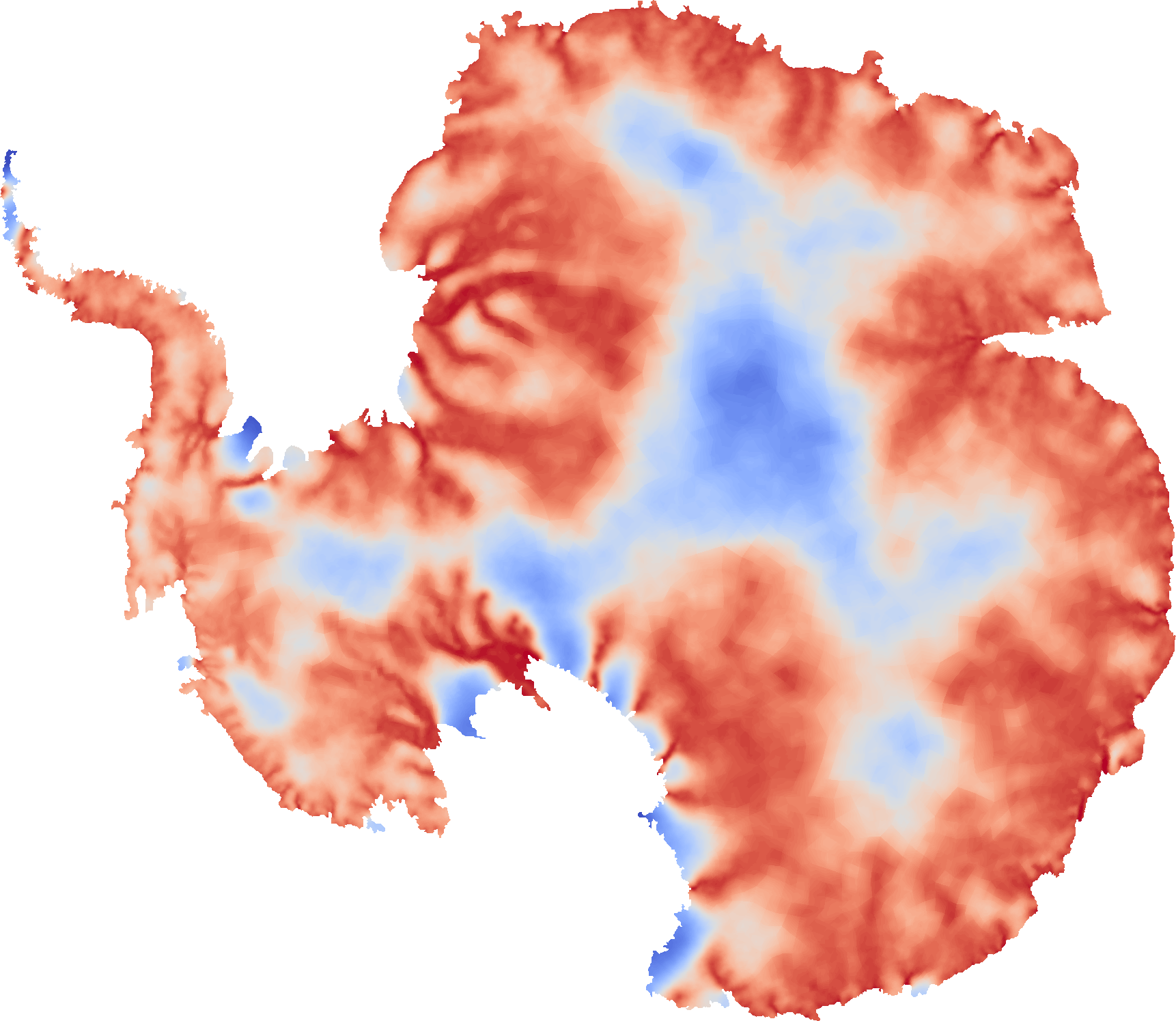}
    \vspace{0.05in} %
    \includegraphics[width=.28\columnwidth]{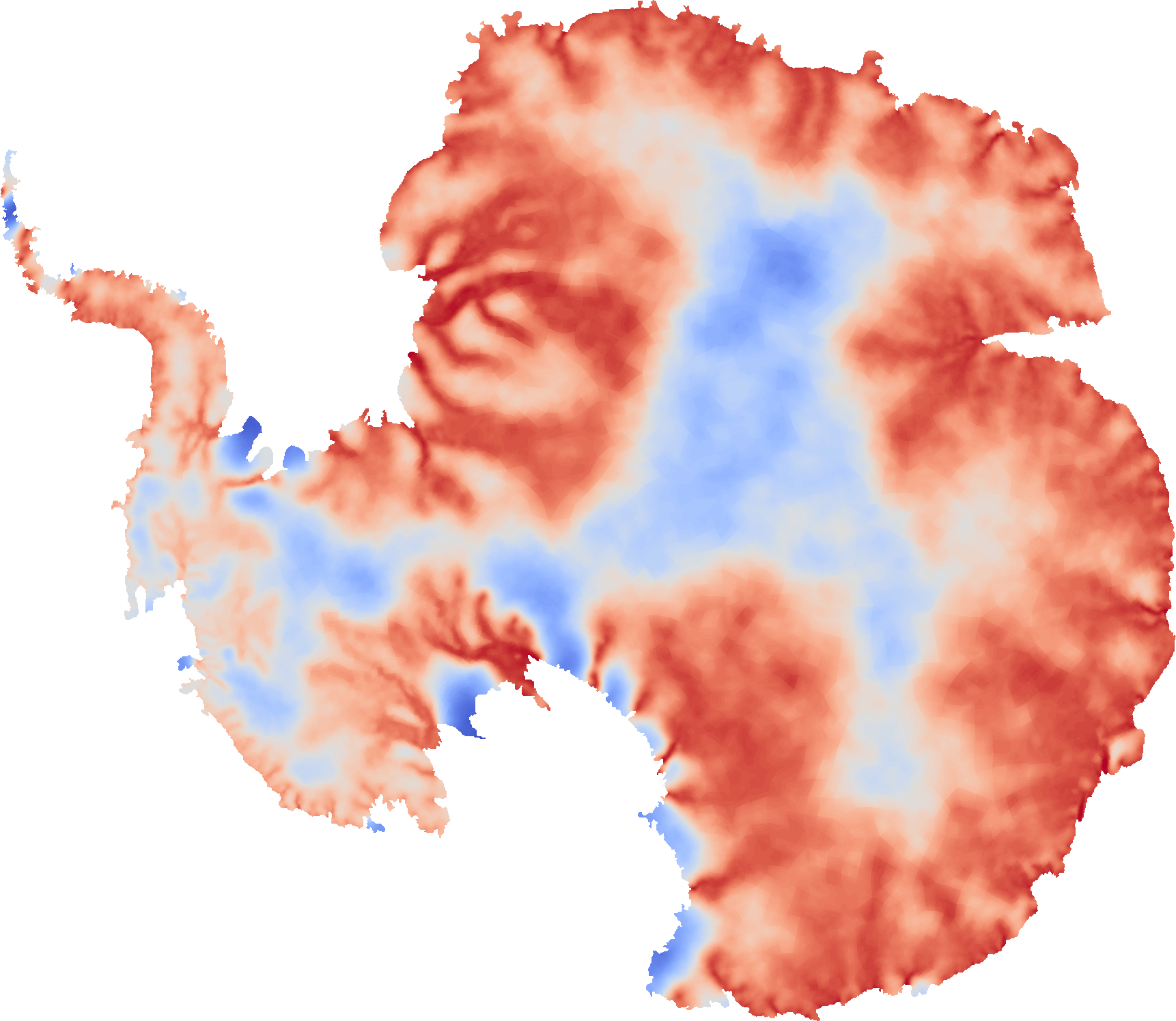} %
  };
  \node (r) at (8.2,-2)
        {\includegraphics[width=0.13\columnwidth]{colorbar_units_samples.png}};
        \draw[fill=white,white] (7,-0.7) rectangle (10,0);
        \node at (8.2,-0.1) {$\exp(\beta)^{1/3}$};
        \node at (8.2,-0.45) {[(Pa/(km/a))$^{1/3}$]};
\end{tikzpicture}
\caption{%
  Samples from the prior (top row) and Gaussian approximation of the
  posterior (bottom row) distributions for the basal sliding parameter
  field.  The difference between the prior and posterior samples
  reflects the information gained from the observational data. The large
  scale features of the posterior samples consistently resemble the
  posterior mean shown in Figure~\ref{fig:mean}.
  Note the small variability in the fast ice flow
  regions, while central and West Antarctica exhibit large variability
  in the inferred basal sliding parameter field.}
\label{fig:samples}
\end{figure}

Figure~\ref{fig:mean} depicts the MAP point,
and Figure~\ref{fig:variance} shows the prior and posterior pointwise
standard deviations. One observes that the uncertainty is vastly reduced
everywhere in the domain, and that the reduction is
greatest along the fast ice flow regions.
In Figure~\ref{fig:samples} we show samples (of the basal sliding
parameter field) from the prior (top row) and from the Gaussian
approximation of the posterior (bottom row) pdf. The difference
between the two sets of samples reflects the information gained from
the data in solving the inverse problem. The differences in the basal
sliding parameter field across the posterior samples demonstrate that
in the fast ice flow regions there is small variability, while in the
center and in West Antarctica, we observe larger variability in the
inferred parameters, reflecting uncertainty due to insensitivity of
surface velocities to the basal sliding in slow velocity and small
velocity gradient regions.

This uncertainty in the inference of the basal sliding parameter
field, however, is merely an intermediate quantity. What is of
ultimate interest is predictions of output quantities of interest,
with associated uncertainties, using the ice sheet model with inferred
parameters and their associated uncertainties, which have been
computed using the methods of this section. This is the subject of the
next section.

\section{Prediction with quantified uncertainty: 
 forward propagation of basal sliding parameter uncertainty to mass flux
 prediction} 
\label{sec:uncertainty_propagation} Once the inverse problem to infer
the unknown basal sliding parameter field from observed surface
velocities has been solved, and the uncertainty in this inference
quantified through a Gaussian approximation of the posterior (made
tractable by a low-rank representation of the prior-preconditioned data
misfit Hessian), we are ready to propagate the uncertain basal sliding
parameter field through the ice flow model to yield a prediction of our
quantity of interest with associated uncertainty.  Ultimately our
interest is in predicting the ice mass flux to the ocean several decades
in the future, under various climate change scenarios.  However, this
requires a model of the ice sheet as an evolving body, more
mechanistic basal boundary conditions, and coupling to atmosphere and
ocean models. 
In the present work, we have instead chosen to illustrate
our data-to-prediction framework with a simpler quantity of interest $Q$
given by the ice mass flux to the ocean using the steady state ice sheet
model employed in the previous sections.  Below, we describe scalable
algorithms for this final step of our data-to-prediction framework.

Formally, this amounts to solving a system of stochastic 
PDEs given by the nonlinear Stokes forward model with the uncertain
basal sliding parameter described by a Gaussian random field. While
the low rank approximation of Section \ref{sec:low-rank} has resulted
in significant dimensionality reduction (from $O(10^6)$ to $O(10^3)$,
as seen in Figure~\ref{fig:spectrum}), the {\em effective} dimension,
$\sim$4000, is still large in absolute terms. Given that this many
modes in parameter space are required to quantify the uncertainty in
the basal sliding coefficient parameters, and given the expense of
solving the large-scale highly-nonlinear forward ice sheet flow
problem, the use of Monte Carlo sampling methods would
be prohibitive, since millions of forward solves would likely be
required to characterize the statistics of the prediction
quantity. Similarly, for the problem we target, the use of polynomial
chaos methods would be prohibitive due to the curse of dimensionality
that afflicts such methods. 

Instead, consistent with our Gaussian approximation of the Bayesian
solution of the inverse problem, and our desire to scale to very large
parameter dimensions, here we linearize the parameter-to-prediction
map at the MAP point, resulting in a Gaussian approximation of the
prediction pdf, $\mathcal{N}(Q_\text{\tiny MAP},\Gamma_{\!
  \text{prediction}})$. The mean of this prediction pdf,
$Q_\text{\tiny MAP}$, is computed by solving the forward ice sheet
flow model \eqref{eq:forward} using $\map$ as the basal sliding
parameter, i.e., the MAP point solution of the inverse problem. The
covariance operator, $\Gamma_{\!  \text{prediction}}$, is found by
propagating the covariance of the model parameters (which is given by
the inverse of the Hessian evaluated at $\map$, i.e.,
$\mathcal{H}^{-1}(\map)$), through the linearized
parameter-to-prediction map, i.e., by 
\begin{equation}
\label{eq:prediction_covariance}
\Gamma_{\! \text{prediction}} := \mathcal{F}(\map) \,
        \mathcal{H}^{-1}(\map) \, \mathcal{F}^*(\map),
\end{equation}
where $\mathcal{F}(\map)$ is the Jacobian of the
parameter-to-prediction map, evaluated at the MAP point $\map$, and
the Hessian at the MAP, $\mathcal{H}(\map)$, is defined
by its action in a direction by \eqref{eq:hessian_action}, which
involves solution of the incremental forward \eqref{eq:incforward} and
incremental adjoint \eqref{eq:incadjoint} Stokes problems.

One of the key ideas to enabling scalability of the
prediction-under-uncertainty problem is that the Jacobian of the
parameter-to-prediction map $\mathcal{F}(\map)$ can be determined for
each prediction quantity $Q$ by computing the gradient of $Q$ with
respect to the parameter field $\beta$. In our case, we are interested
in using the steady state ice flow model to predict the net ice mass
flux into the ice shelves, and eventually into the ocean,
\begin{equation}
\label{qoi}
Q(\beta) := \int_{\Gamma_o} \!\! \rho \bs{u}(\beta) \cdot \bs{n} \, ds,
\end{equation}
where $\Gamma_o$ is an outflow boundary of interest. 
The gradient of $Q$ with respect to $\beta$ evaluated at $\map$ can
then be found as follows. First, solve the forward problem
\eqref{eq:forward} with basal sliding parameter field given by $\map$. Then,
solve an adjoint problem defined for the quantity $Q$, i.e.,
\begin{subequations}\label{eq:adjoint_prediction}
\begin{alignat}{2}
  \gbf{\nabla} \cdot \vv &= 0  &\; &  \quad \text{ in }  \Omega\\
  - \gbf{\nabla} \cdot \gbf{\sigma}_{\!\vv} &= \gbf{0} & & \quad
    \text{ in } \Omega\\ 
 \gbf{\sigma}_{\!\vv} \gbf{n} &= 0  & \; & \quad \text{ on }
 \Gamma_{\!\mbox{\tiny} t} \\ 
 \gbf{\sigma}_{\!\vv} \gbf{n} &= - \rho \gbf{n}  & \; & \quad \text{ on }
 \Gamma_{\!\mbox{\tiny} o} \\ 
  T\gbf{\sigma}_{\!\vv} \gbf{n} + \exp(\map) \gbf T\vv &=
  \gbf{0}, \;\; \vv\cdot \gbf{n} = 0 &\; & \quad \text{ on }
  \Gamma_{\!\mbox{\tiny} b}
\end{alignat}
\end{subequations}
where the adjoint stress is given by~\eqref{eq:tensor-viscosity}.
Since we are using the same ice flow model for the prediction problem
that was used to infer $\beta$, the adjoint problem
\eqref{eq:adjoint_prediction} for $Q$ resembles the adjoint problem
for the regularized data misfit functional $\mathcal{J}$ in
\eqref{objfunction}, but with a different source term given by the
variation of $Q$ with respect to the state variables
$(\gbf{u},p)$. This means that the adjoint Stokes solver can be reused
for this purpose (which in turn is the same as the linearized forward Stokes
solver). 
Once the forward velocity $\gbf{u}$ and adjoint velocity $\gbf{v}$ are found,
$\mathcal{F}(\map)$, which is the gradient of $Q$ at $\map$, can be
found from the gradient expression \eqref{eq:gradient} without the
regularization term, i.e., from 
\begin{equation}\label{eq:predgrad}
  \mathcal{F}(\map):= \exp(\map)\; \gbf T \uu \cdot
  \gbf T \vv  \quad \text{ on
  } \Gamma_{\smallsub b}.
\end{equation}
Note that if a different flow model is used for the prediction phase
(such as one governing a dynamically evolving ice sheet), then a new
adjoint equation and gradient expression will have to be derived for
that model for the prediction step.

Now that $\mathcal{F}(\map)$ has been found, the next step is to form
$\mathcal{H}^{-1}(\map) \mathcal{F}^*(\map)$, which could be found by solving
a linear system using the preconditioned CG method described in
Section \ref{ipmethod}. This would require a number of
forward/adjoint incremental Stokes solves equal to the number of CG
iterations. Instead, since $\mathcal{H}^{-1}(\map)$ is available in the
compact form given by \eqref{eqn:reduced_post} (based on the low-rank
approximation of the prior-preconditioned Hessian of the data misfit),
the product $\mathcal{H}^{-1}(\map) \mathcal{F}^*(\map)$ can be formed without
having to solve incremental Stokes equations, at the cost of just
linear algebra (vector scalings, additions and inner products). The final step of the
inner product of $\mathcal{F}(\map)$ with $\mathcal{H}^{-1}(\map)
\mathcal{F}^*(\map)$ is straightforward. The result is the covariance of
the ice mass flux $Q$, given uncertainties in the satellite
observations of surface velocity, uncertainties in the basal sliding
parameter field as inferred from the satellite data, and the ice sheet
flow model and corresponding prediction quantity $Q$. \note{For a single
$Q$}, the prediction density $\mathcal{N}(Q_\text{\tiny MAP},\Gamma_{\!
  \text{prediction}})$ is univariate. The cost to obtain this pdf,
once the low rank-based inverse Hessian \eqref{eqn:reduced_post} has
been computed in the Bayesian inference phase, is just a forward solve
followed by an (as always, linear) adjoint solve, i.e., at little cost
beyond the forward solve.

\begin{figure}[ht]
  \centering
\begin{minipage}{1\textwidth}
  \includegraphics[width=.328\columnwidth]{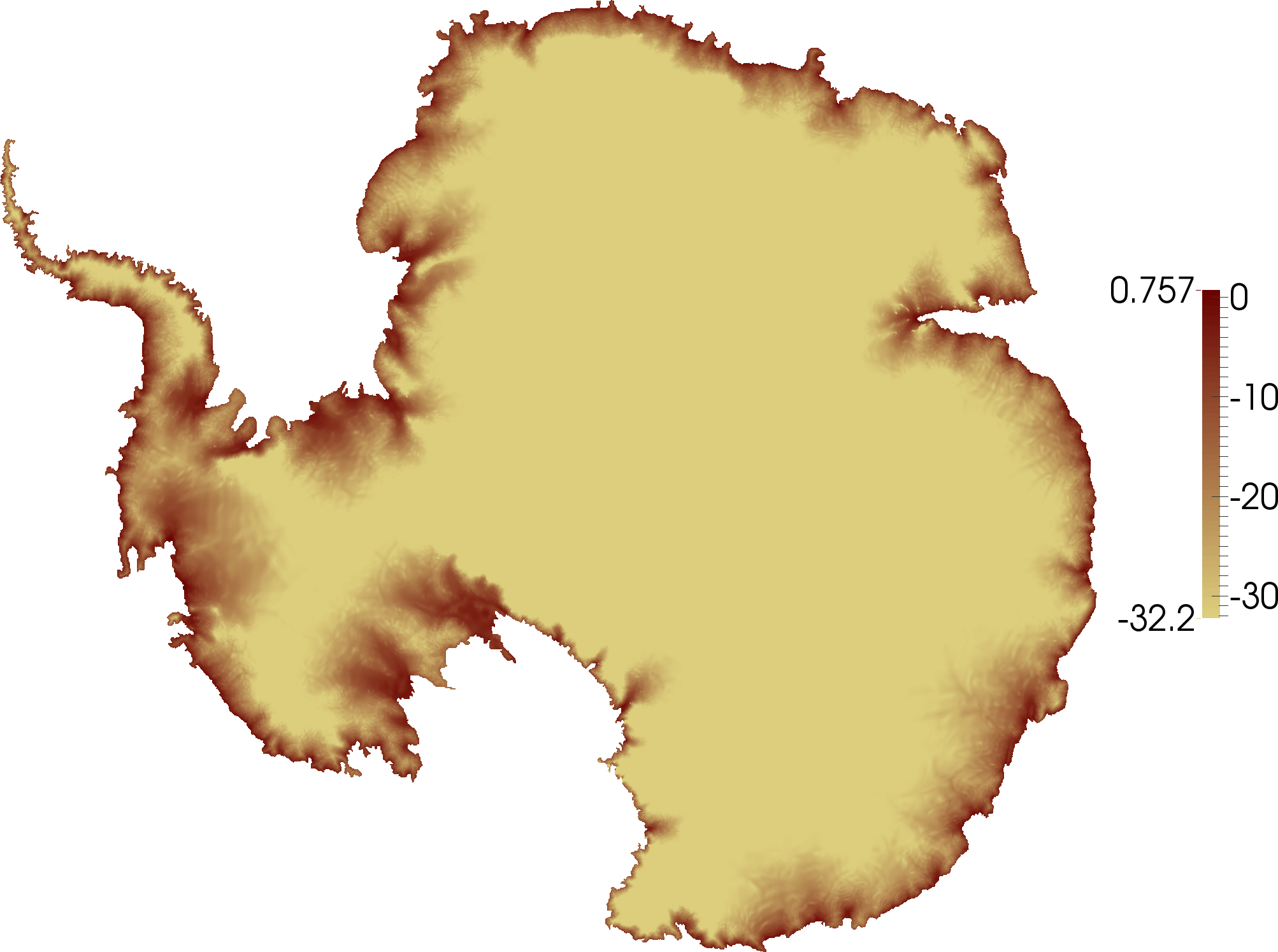} \hfill
  \includegraphics[width=.328\columnwidth]{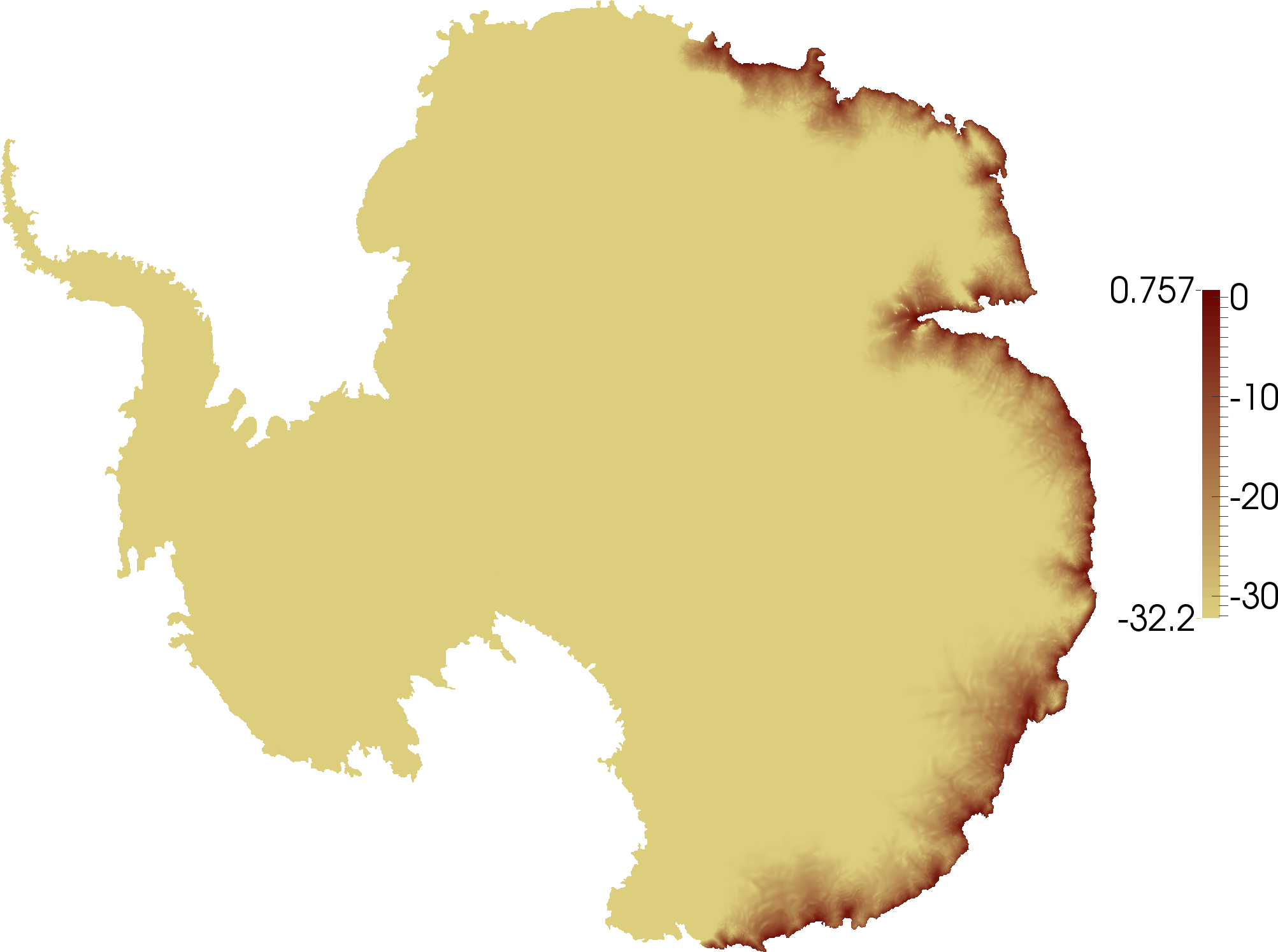} \hfill
  \includegraphics[width=.325\columnwidth]{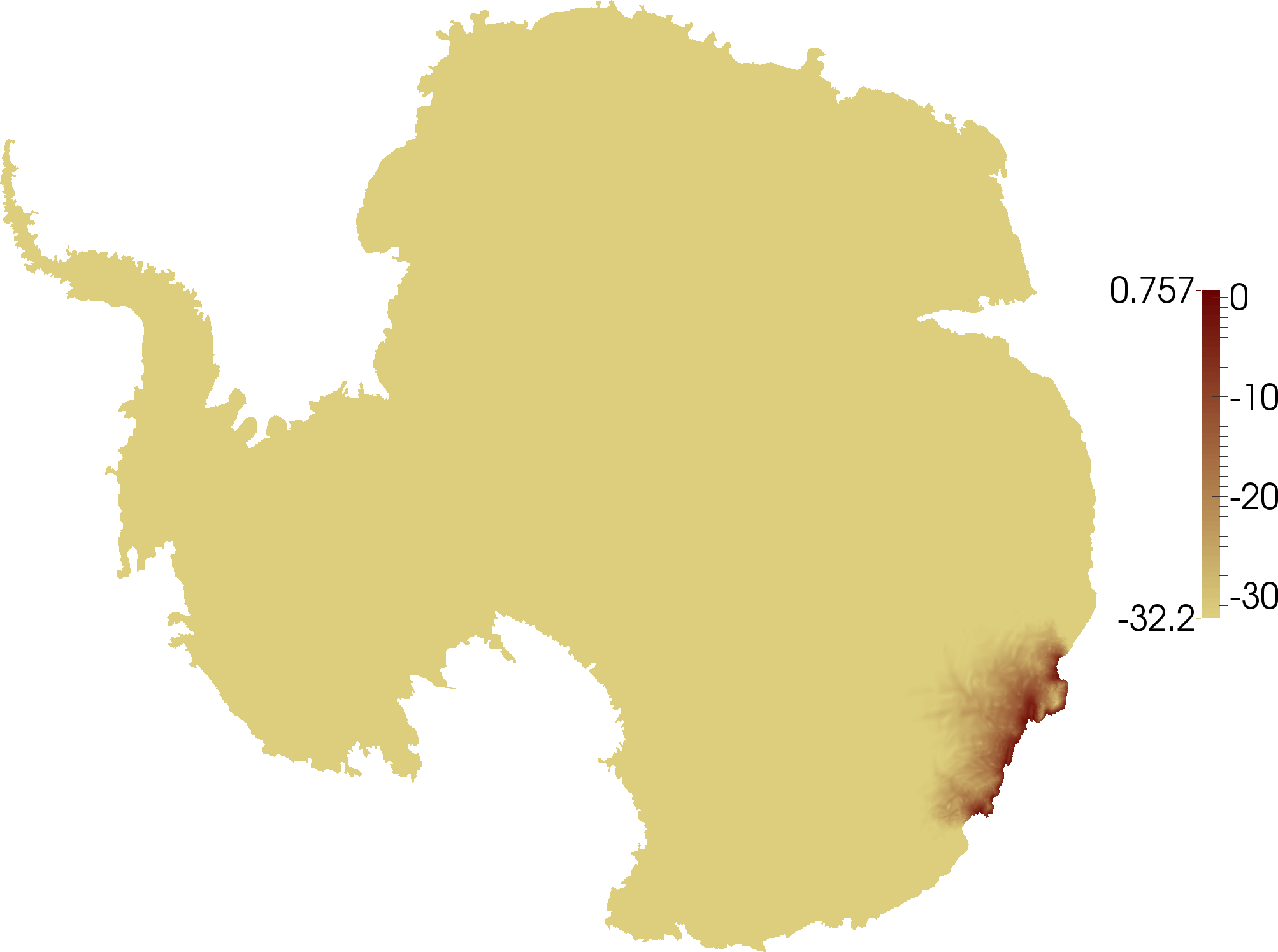}\\
  \includegraphics[width=.34\columnwidth]{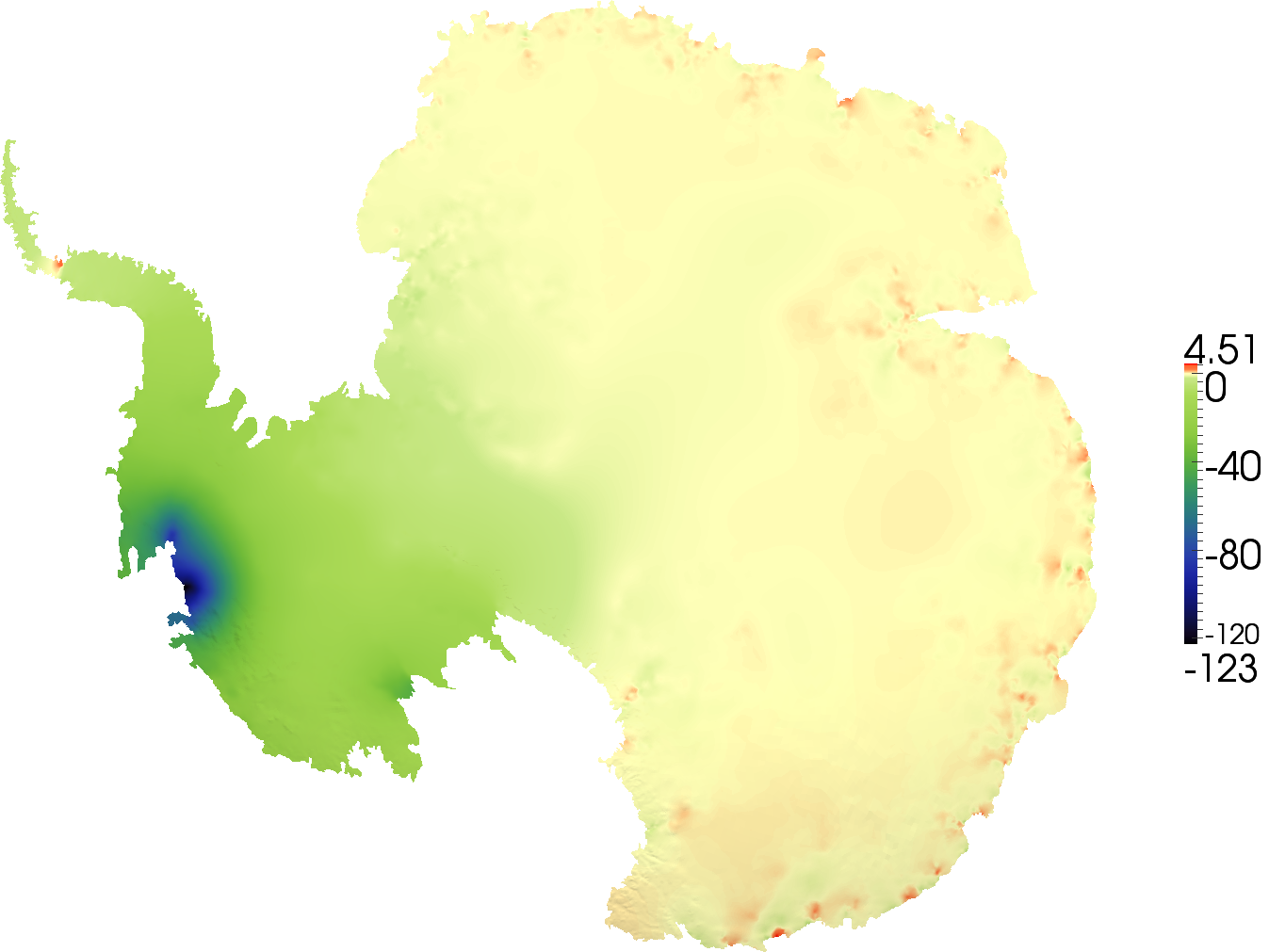} \hfill
  \includegraphics[width=.33\columnwidth]{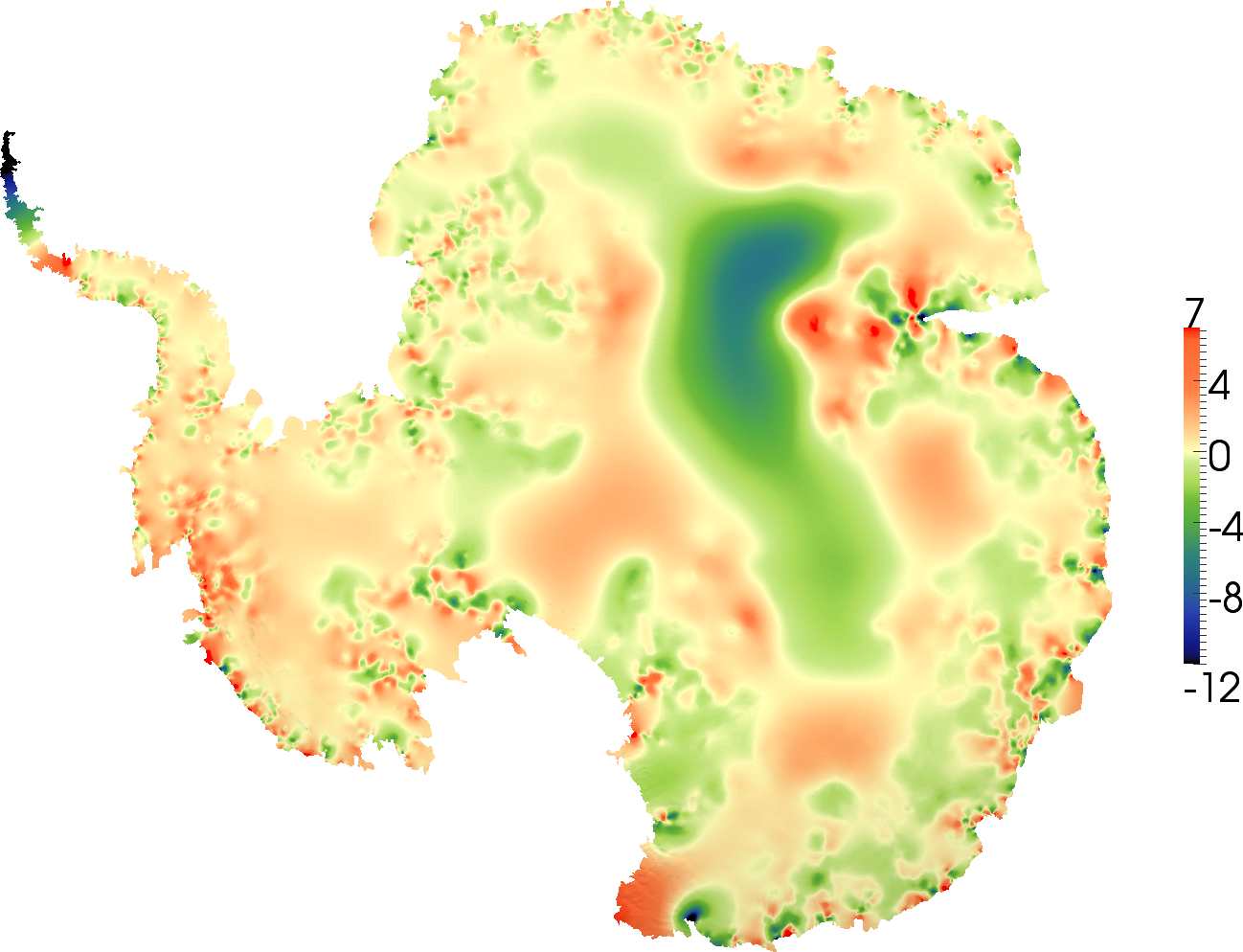}\hfill
  \includegraphics[width=.31\columnwidth]{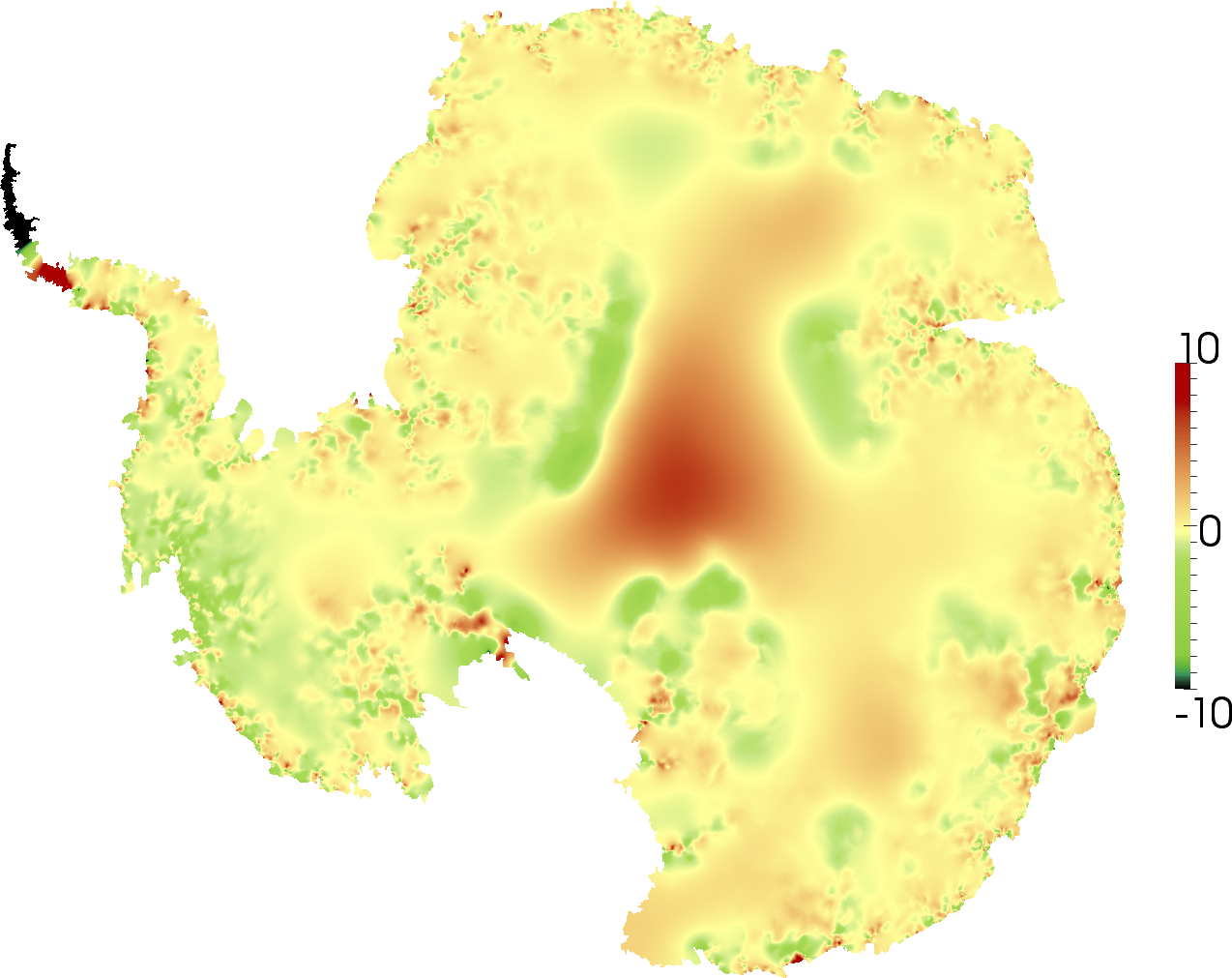}  \hfill
\end{minipage}
\caption{The (logarithm of the) gradient of the prediction quantity
  $Q$ with respect to
  the uncertain basal sliding parameter $\beta$ (top row) and the
  directions that jointly maximize uncertainty and sensitivity for the
  prediction (bottom row) following \cite{LiebermanWillcox12}. For
  visualization purposes, in the top row we plot $\ln
  (|\mathcal{F}(\map)| + 10^{-10})$. The quantities of interest are
  the mass flux of ice into the ocean for the entire boundary (left
  column), for all of East Antarctica (center column), and for the
  Totten Glacier region (right column). The mean and standard
  deviation of the prediction probability distribution for the three
  ice mass fluxes are: 1170.83 $\pm$ 1080.09, 359.60 $\pm$ 1.02, and
  71.24 $\pm$ 0.30 Gt/a, respectively.}

\label{fig:pred}
\end{figure}

Figure~\ref{fig:pred} provides results of uncertainty propagation from
parameters to predictions following the framework outlined above.
Recent work on ``inference for prediction'' has shown that, for linear
inverse problems and linear parameter-to-prediction maps, there is a
unique direction in parameter space that influences the prediction
quantity of interest \cite{LiebermanWillcox12}.
Finding this direction involves identifying parameter modes that are
both informed by the observational data and also required for
estimating the quantity of interest. We adopt the
inference-for-prediction (IFP) algorithm presented more generally for
multiple quantities of interest
in~\cite{LiebermanWillcox12}. Denoting the %
Hessian square root by $\mathcal{G} := \mathcal{H}^{-1/2}(\map)$ and using
the linearized parameter-to-prediction map $\mathcal{F}(\map)$,
we compute the
eigendecomposition $\PP \SI^2 \PP^*$ of ${\mathcal G}^* {\mathcal
  F^*(\map)} \mathcal F(\map) \mathcal G$.
Since this
operator has rank 1, the eigenvector $\PP$ corresponding
to the only nonzero eigenvalue is given by ${\mathcal G}^* {\mathcal{F}^*(\map)}$, and the
corresponding eigenvalue is
$\SI^2 = \|{\mathcal G}^*
{\mathcal F^*(\map)}\|^2$.
Following \cite{LiebermanWillcox12}, the
influential direction for prediction (based on linearizations of the
parameter-to-observable map and the parameter-to-prediction map) is 
given by
\begin{align*}
  \WW = \mathcal G \PP \SI^{-1/2}.
\end{align*}
Using the explicit form of $\PP$, the influential direction for
prediction in the IFP algorithm thus simplifies to 
$\WW = \mathcal G {\mathcal G}^* {\mathcal{F}^*(\map)} \SI^{-1/2} =
  \SI^{-1/2} \mathcal H(\map)^{-1}\mathcal F^*(\map)$.
This direction, or ``mode'' of the parameter field, is depicted in the
bottom row of Figure~\ref{fig:pred} for ice mass fluxes over three
different portions of the outflow boundary (all of Antarctica, just
East Antarctica, and just the Totten Glacier).
The mean and standard deviation of the prediction probability
distribution for the three ice mass fluxes are: 1170.83 $\pm$ 1080.09,
359.60 $\pm$ 1.02, and 71.24 $\pm$ 0.30 Gt/a, respectively. We also
computed standard deviations with the prior covariance and, as expected,
found larger values, namely 1327.22, 525.24, and 179.73 Gt/a for all
of Antarctica, East Antarctica and the Totten Glacier, respectively.
The top row of Figure~\ref{fig:pred}
portrays the gradient of $Q$ with respect to $\beta$, computed
using~\eqref{eq:predgrad} for each case; these plots show the regions
in parameter space to which $Q$ is most sensitive.  
Note the differences with the bottom row, which captures not only
sensitivity of the prediction quantity to $\beta$ (i.e., the gradient
$\mathcal{F}$), but also uncertainty in the $\beta$ field (i.e.,
$\mathcal{H}^{-1}$).
The most uncertain regions are not necessarily the most sensitive, and
vice versa. 
Since, as implied by Figure~\ref{fig:samples}, central and West
Antarctica 
exhibit the largest uncertainties, we expect these regions to play an
important role in the parameter modes shown in the bottom row of
Figure~\ref{fig:pred}, especially where the sensitivity with respect to the
basal sliding parameter is not dominant, e.g., in the bottom left figure. 
On the other hand, in the center and right figures in the bottom row,
focusing on East Antarctica, where the posterior
samples suggest lower uncertainties, the modes show mixed uncertainty
and sensitivity influence.

\section{Conclusions}
\label{sec:conclusions}

We have presented a scalable framework for solving end-to-end,
data-to-prediction problems, presented in the context of the flow of
the Antarctic ice sheet, and motivated by prediction of its
contribution to sea level.  We begin with observational data and
associated uncertainty, then infer model parameters from the data
(Section \ref{sec:inverse}), quantify uncertainties in the inference
of the parameters (Section \ref{sec:bayesian}), and finally propagate
those uncertain model parameters to yield a prediction quantity of
interest with quantified uncertainties (Section
\ref{sec:uncertainty_propagation}). {\em We show that the cost of the
  entire data-to-prediction framework, when measured in forward or
  adjoint Stokes solves, is a constant independent of the parameter or
  data dimension.} When combined with a scalable forward solver such
as that presented in Section \ref{sec:forward}, this results in a
data-to-prediction framework that is independent of the state
dimension and number of processor cores as well.

This scalability is a consequence of three properties of the
data-to-prediction process and their exploitation by our algorithms:
(1) when inferring the parameter field, the data are informative about
only a low-dimensional subspace within the high-dimensional parameter
space, and thus a Newton-CG optimization method, preconditioned by the
regularization operator, converges in a number of Newton and CG
iterations that is independent of the parameter and data dimensions,
depending only on the information content of the data; (2) when
estimating the uncertainty in the inverse solution, the same property
dictates that the Hessian of the data misfit admits a low rank
representation, and this can be extracted via randomized SVD in a
number of matrix-free Hessian-vector products (each of which requires
a pair of incremental forward/adjoint Stokes solves) that is also
independent of the data and parameter dimensions, and again depends
only on the information contained within the data; and (3) when
propagating the inferred parameter uncertainties forward through the
ice flow model to yield predictions with quantified uncertainties, the
prediction pdf can be formed through the action of the inverse Hessian
on the gradient of the prediction quantity with respect to the
uncertain parameters, which in turn is found through an additional
adjoint Stokes solve. Thus, {\em the entire data-to-prediction process
  is sensitive only to the true information contained within the data,
  as opposed to the ostensible data or parameter dimensions.}

The uncertainty analysis presented here relies on linearizations of
the parameter-to-observable and parameter-to-prediction maps, leading
to Gaussian approximations of the parameter posterior pdf and
prediction quantity of interest pdf. Ultimately one would like to
relax these approximations and fully explore the resulting
non-Gaussian pdfs; how to do this while retaining scalability for such
large-scale complex problems remains an open question, and is a
subject of ongoing work.

\section*{Acknowledgment}
We thank the anonymous reviewers for their thorough reading of
the manuscript, and for their valuable comments and suggestions
that improved the quality of the paper. Support for this work was
provided by: the U.S. Air Force Office of Scientific Research
Computational Mathematics program under award number FA9550-12-1-0484;
the U.S. Department of Energy Office of Science Advanced Scientific
Computing Research program under award numbers DE-FG02-09ER25914,
DE-FC02-13ER26128, and DE-SC0010518; and the U.S. National Science
Foundation Cyber-Enabled Discovery and Innovation program under
awards CMS-1028889 and OPP-0941678. This research used resources
of the Oak Ridge Leadership Facility at the ORNL, which is
supported by the Office of Science of the DOE under Contract
No. DE-AC05-00OR22725. Computing time on TACC's Stampede was
provided under XSEDE, TG-DPP130002.
\section*{References}

\end{document}